%% file: main.tex
\documentclass[journal,twoside]{IEEEtran}
\IEEEoverridecommandlockouts
\usepackage{cite}

\ifCLASSINFOpdf
  \usepackage[pdftex]{graphicx}
 \graphicspath{{images/}}
  \DeclareGraphicsExtensions{.pdf,.jpeg,.png}
\else
  \usepackage[dvips]{graphicx}
  \graphicspath{{images/}}
  \DeclareGraphicsExtensions{.eps}
\fi
\usepackage[cmex10]{amsmath}
\usepackage{amssymb}

\ifCLASSOPTIONcompsoc
  \usepackage[caption=false,font=normalsize,labelfont=sf,textfont=sf]{subfig}
\else
  \usepackage[caption=false,font=footnotesize]{subfig}
\fi

\usepackage{abbrevs}
\usepackage{color}
\usepackage{comment}


\begin{document}
\title{Likelihood Functions with Parameter-Dependent Support: A Survey of the Cram\'{e}r-Rao-Leibniz Lower Bound}
\author{
\IEEEauthorblockN{Qin~Lu$^\dag$, Yaakov Bar-Shalom$^\ddag$, Peter Willett$^\ddag$\\}
\IEEEauthorblockA{$^\dag$School of Electrical and Computer Engineering, University of Georgia, USA\\}
\IEEEauthorblockA{$^\ddag$Department of Electrical and Computer Engineering, University of Connecticut, USA\\}
\thanks{Q. Lu was supported by NSF CAREER \#2340049.}
}

\maketitle
\begin{abstract}
Parameter estimation is a fundamental problem in science and engineering. In many safety-critical applications, one is not only interested in a {\it point} estimator, but also the uncertainty bound that can self-assess the accuracy of the estimator. In this regard, the Cram\'{e}r-Rao lower bound (CRLB) is of great importance, as it provides a lower bound on the variance of {\it any} unbiased estimator. In many cases, it is the only way of evaluating, without recourse to simulations, the expected accuracy of numerically obtainable estimates. For the existence of the CRLB, there have been widely accepted regularity conditions, one of which is that the support of the likelihood function (LF) -- the pdf of the observations conditioned on the parameter of interest-- should be independent of the parameter to be estimated. This paper starts from reviewing the derivations of the classical CRLB under the condition that the LF has parameter-independent support. To cope with the case of parameter-dependent support, we generalize the CRLB to the {\it Cram\'{e}r-Rao-Leibniz lower bound (CRLLB)}, by leveraging the general Leibniz integral rule. Notably, the existing results on CRLLB and CRLB are unified under the framework of CRLLB with multidimensional parameters. Then, we survey existing examples of LFs to illustrate the usefulness of the CRLLB in providing valid covariance bounds.




\end{abstract}


\section{Introduction}

In estimation theory and statistics, the Cram\'{e}r-Rao lower bound (CRLB), as a lower bound on the variance of {\it any} unbiased estimator, can provide a bound to self-assess the estimation accuracy, without recourse to simulations~\cite{BarShalomLK:2001, VanTrees:2013}.  In several nonlinear problems, where the ML estimates could be obtained only numerically  \cite{osborne2013statistical, lu2017tracking, lu2018efficient, yang2019estimation}, it was shown via hypothesis testing that the CRLB could be accepted as the actual covariance of the corresponding estimator even for a small number of samples, i.e., they were efficient. Further, it provides an important performance measure to guide the design of system parameters for various applications in signal processing~\cite{liu2021cramer, lu2018measurement, lu2019measurement}. 


To estimate the unknown parameter $\mathbf{x}$ from the measurement $\mathbf{z}$, we will rely on the likelihood function (LF) $p(\mathbf{z}|\mathbf{x})$, which, as a function of $\mathbf{x}$, represents the pdf of $\mathbf{z}$ conditioned on $\mathbf{x}$.
In order for the classical CRLB to hold, there are the following two requirements\cite{BarShalomLK:2001} \cite{BarShalomWT:2011} \cite{VanTrees:2013}:
\begin{itemize}
\item Existence and absolute integrability with respect to (w.r.t.) $\mathbf{z}$ of the first two derivatives of the LF $p(\mathbf{z}|\mathbf{x})$ w.r.t. $\mathbf{x}$; and
\item The support of the LF $p(\mathbf{z}|\mathbf{x})$ be independent of $\mathbf{x}$.
\end{itemize}

For the case of parameter-dependent LF support, which arises when an unknown parameter is observed in the presence of additive noise with a finite support pdf, it was until recently believed that one can make no informed decision concerning the existence of the CRLB. Under careful scrutiny, however, it was stated in \cite{YBS_2014_CDC} \cite{YBS_2014_July} \cite{YBS_2014_March} that the above regularity conditions are too stringent: the following two relaxed conditions were presented instead as necessary and sufficient for the CRLB to hold:
\begin{itemize}
\item The expected value of the square of the first partial derivative of the log-LF (LLF) $\ln p(\mathbf{z}|\mathbf{x})$ w.r.t. $\mathbf{x}$, conditioned on $\mathbf{x}$, is finite; and
\item The support of the LF $p(\mathbf{z}|\mathbf{x})$ is independent of the parameter $\mathbf{x}$; or if it depends on $\mathbf{x}$, it has to be continuous at the boundary of its support (equals to zero at the boundary of its support).
\end{itemize}
The former does not require the twice differentiability of the LF as in the previous literature, but only almost surely (a.s.) the first differentiability of LLF (a.s. is w.r.t. the probability measure given by the LF). An interesting case in point for this is the Laplace LF which, while not differentiable at the origin, is a.s. differentiable because the point where it is not differentiable has probability zero and thus it does not contribute to the (Lebesgue) integral that yields the FIM. Also, note that the Laplace LF belongs to the exponential family but still it does not have an efficient estimator --- this is because belonging to the exponential family is only a necessary but not a sufficient condition for an efficient estimator to exist.

For the case where the LF is not continuous at the end points of its support, the authors in \cite{YBS_2014_CDC} \cite{YBS_2014_July} \cite{YBS_2014_March} derived a new modified bound --- designated as the Cram\'{e}r-Rao-Leibniz lower bound (CRLLB) since it relies on the Leibniz integral rule. The CRLLB provides, among other results, a valid bound for a uniformly distributed measurement noise, which has been a long-standing problem because the CRLB was believed not to hold for it. In the literature, variants of Cram\'{e}r-Rao type lower bound have been investigated, for example, a Chapman-Robbins form of the Barankin bound was used to derive a multiparameter CRLB on the estimation error covariance when the parameter to be estimated was constrained to lie in a subset of the parameter space in \cite{gorman1990lower}.
 To the best of the authors' knowledge, there is no work in the literature dealing with LF with finite support, except for \cite{YBS_2014_July} \cite{YBS_2014_CDC}\cite{YBS_2014_March}. However, these references were restricted to scalar parameters. 
 Later, the CRLLB for multidimensional parameters is obtained by applying the general Leibniz integral rule using the ``interior product"~\cite{lu2017multidimensional, lu2017multidimensional_fusion}.
It is worth mentioning that the scalar CRLLB is based on the Cauchy-Schwarz inequality for a scalar-valued function, while the {\em multidimensional} CRLLB derivation relies on the positive definiteness of the diagonal blocks of the partition of the inverse of a partitioned positive definite matrix. The scalar CRLLB of an estimator is attained (i.e., the estimator is efficient) if it meets the Cauchy-Schwarz collinearity condition, while for the {\em multidimensional} CRLLB case, the regularity condition for an efficient estimator necessitates the {\it generalized collinearity condition}. 

 This paper surveys existing results on CRLB and CRLLB, which are unified under the framework of CRLLB for multidimensional parameters. Specifically, the rest of the paper is organized as follows. Section \ref{sec:CRLB} reviews the derivation of the classical CRLLB with parameter-independent support. To account for parameter-dependent support, Section \ref{sec:CRLLB} presents the CRLLB for multidimensional parameters, which, as shown later, unifies existing results on CRLB and CRLLB. Sections \ref{sec:examples}-\ref{sec:uniform_example} provide some illustrative examples of LFs with parameter-dependent support for different dimensional measurements. Concluding remarks and discussions are given in Section \ref{sec:conclusion}. 

\noindent {\bf Notational conventions.} Bold lowercase letters are used to represent vectors with its entries denoted by lowercase plain counterpart with corresponding indices. Matrices are represented by bold capital letters.

\section{Review of the classical CRLB}\label{sec:CRLB}
Given the $n_{\bf z}$-dimensional observation vector ${\bf z}\in \mathbb{R}^{n_{\bf z}}$, our goal is to estimate the $n_{\bf x}$-dimensional parameter vector ${\bf x}$. The connection between ${\bf x}$ and ${\bf z}$ is captured by the so-termed LF $p(${\bf z}$|${\bf x}$)$, which represents the pdf of ${\bf z}$ conditioned on ${\bf x}$. Let $\hat{\mathbf{x}}(\mathbf{z})$ denote an unbiased estimator of ${\bf x}$. Then, the unbiasedness condition yields
\begin{align}
    \mathbb{E}[(\hat{\mathbf{x}}(\mathbf{z}) - \mathbf{x})|\mathbf{x}] = \int_{\mathbb{R}^{n_{\bf z}}}p(\mathbf{z}|\mathbf{x})(\hat{\mathbf{x}}(\mathbf{z}) - \mathbf{x}) d\mathbf{z} ={\bf 0}_{n_{\bf x}\times 1} \label{eq:unbiased_CRLB}
\end{align}
where ${\bf 0}_{n_{\bf x}\times 1}$ is an $n_{\bf x}\times 1$ all-zero vector.

Taking the gradient of ~\eqref{eq:unbiased_CRLB} with respect to (w.r.t.) ${\bf x}$ yields
\begin{align}
& \nabla_{\bf x} \mathbb{E}[(\hat{\mathbf{x}}(\mathbf{z}) - \mathbf{x})|\mathbf{x}] \label{eq:gradient}\\
 & =  \int_{\mathbb{R}^{n_{\bf z}}} \!\!\left(\nabla_{\bf x}p(\mathbf{z}|\mathbf{x})(\hat{\mathbf{x}}(\mathbf{z}) - \mathbf{x})^\top\! \!- p(\mathbf{z}|\mathbf{x}){\bf I}_{n_{\bf x}} \right) d\mathbf{z} = {\bf 0}_{n_{\bf x}} \nonumber
\end{align}
where ${\bf I}_{n_{\bf x}}$ is an $n_{\bf x}\times n_{\bf x}$ identity matrix, ${\bf 0}_{n_{\bf x}}$ is an $n_{\bf x}\times n_{\bf x}$ all-zero matrix, and $\top$ represents the transpose of a vector or matrix. Here, the gradient comes from the integrand {\it only} as the integration bounds in~\eqref{eq:unbiased_CRLB} do {\it not} depend on ${\bf x}$.

Since $\int_{\mathbb{R}^{n_{\bf z}}}p(\mathbf{z}|\mathbf{x}) d\mathbf{z} = 1$, the following equality holds
\begin{align}
   \int_{\mathbb{R}^{n_{\bf z}}} (\nabla_{\bf x}p(\mathbf{z}|\mathbf{x}))(\hat{\mathbf{x}}(\mathbf{z}) - \mathbf{x})^\top d\mathbf{z} = {\bf I}_{n_{\bf x}} \label{eq:gra_p}
\end{align}
Using
\begin{align}
  \nabla_{\bf x}p(\mathbf{z}|\mathbf{x}) =  (\nabla_{\bf x}\ln p(\mathbf{z}|\mathbf{x})) p(\mathbf{z}|\mathbf{x})\label{eq:equality}
\end{align}
it follows that
\begin{align}
 \int_{\mathbb{R}^{n_{\bf z}}} (\nabla_{\bf x}\ln p(\mathbf{z}|\mathbf{x})) p(\mathbf{z}|\mathbf{x})(\hat{\mathbf{x}}(\mathbf{z}) - \mathbf{x})^\top d\mathbf{z} = {\bf I}_{n_{\bf x}}
\end{align}
Define
\begin{subequations}
\begin{align}
\boldsymbol{\epsilon}(\mathbf{z},\mathbf{x}) &\triangleq \hat{\mathbf{x}}(\mathbf{z}) - \mathbf{x} \label{eq: epsilon}\\
\boldsymbol{\gamma}(\mathbf{z},\mathbf{x}) &\triangleq \nabla_\mathbf{x} \ln \,p(\mathbf{z}|\mathbf{x})\label{eq: gamma}
\end{align}
\end{subequations}
Leveraging the Cauchy-Schwarz inequality, one can obtain that
\begin{align}
&\int_{\mathbb{R}^{n_{\bf z}}}  \boldsymbol{\epsilon}(\mathbf{z},\mathbf{x}) \boldsymbol{\epsilon}(\mathbf{z},\mathbf{x})^\top p(\mathbf{z}|\mathbf{x}) d\mathbf{z} \int_{\mathbb{R}^{n_{\bf z}}}  \boldsymbol{\gamma}(\mathbf{z},\mathbf{x}) \boldsymbol{\gamma}(\mathbf{z},\mathbf{x})^\top p(\mathbf{z}|\mathbf{x}) d\mathbf{z} \nonumber\\
&\geq  \left(\int_{\mathbb{R}^{n_{\bf z}}} \boldsymbol{\epsilon}(\mathbf{z},\mathbf{x}) \boldsymbol{\gamma}(\mathbf{z},\mathbf{x})^\top p(\mathbf{z}|\mathbf{x}) d\mathbf{z}\right)^2 =  \mathbf{I}_{n_{\bf x}}
\end{align}
With the estimator's covariance matrix $\mathbf{P}(\mathbf{x})$ and the Fisher information matrix (FIM) $\mathbf{J}(\mathbf{x})$ being denoted as
\begin{align}
\mathbf{P}(\mathbf{x}) &\triangleq 
 \mathbb{E}[\boldsymbol{\epsilon}(\mathbf{z},\mathbf{x})\boldsymbol{\epsilon}(\mathbf{z},\mathbf{x})^{\top}|\mathbf{x}]\label{eq:P} \\
\mathbf{J}(\mathbf{x}) &\triangleq 
 \mathbb{E}[\boldsymbol{\gamma}(\mathbf{z},\mathbf{x})\boldsymbol{\gamma}(\mathbf{z},\mathbf{x})^{\top}|\mathbf{x}] \label{eq:J}
\end{align}
it is apparent that 
\begin{align}
\mathbf{P}(\mathbf{x}) \geq \mathbf{J}(\mathbf{x})^{-1}
\end{align}
where the equality holds if and only if the generalized collinearity condition holds
\begin{align}
\ \hat{\mathbf{x}}(\mathbf{z}) - \mathbf{x} = {\bf J}(\mathbf{x})^{-1}\nabla_\mathbf{x} \ln \,p(\mathbf{z}|\mathbf{x})\qquad  a.s.\forall \mathbf{z} \in \mathbb{R}
\end{align}
The inverse of FIM $\mathbf{J}(\mathbf{x})^{-1}$ is the well-known CRLB for the covariance of {\it any} unbiased estimators.

Next, we will derive an alternative form of the FIM and CRLB.
With $\int_{\mathbb{R}^{n_{\bf z}}} p(\mathbf{z}|\mathbf{x}) d\mathbf{z} = 1$, it holds that
\begin{align}
   \int_{\mathbb{R}^{n_{\bf z}}}\!\!\!\nabla_{\bf x} p(\mathbf{z}|\mathbf{x}) d\mathbf{z} =  \int_{\mathbb{R}^{n_{\bf z}}} \!\! p(\mathbf{z}|\mathbf{x}) \nabla_{\bf x} \ln p(\mathbf{z}|\mathbf{x}) d\mathbf{z} ={\bf 0}_{n_{\bf x}\times 1}
\end{align}
Further, taking the gradient w.r.t. ${\bf x}$ and using~\eqref{eq:equality} yield
\begin{align}
& \int_{\mathbb{R}^{n_{\bf z}}} p(\mathbf{z}|\mathbf{x}) \nabla_{\bf x}\nabla_{\bf x}^\top \ln p(\mathbf{z}|\mathbf{x})d\mathbf{z} \nonumber\\
& =- \int_{\mathbb{R}^{n_{\bf z}}}p(\mathbf{z}|\mathbf{x})\nabla_{\bf x} \ln p(\mathbf{z}|\mathbf{x}) \nabla_{\bf x}^\top \ln p(\mathbf{z}|\mathbf{x}) d\mathbf{z} \label{eq:FIM_hessian}
\end{align}
thus obtaining the alternative form of the FIM as the expectation of the negative Hessian of the LLF.

To sum up, the existence of classical CRLB hinges on the following two requirements\cite{BarShalomLK:2001} \cite{BarShalomWT:2011} \cite{VanTrees:2013}:
\begin{enumerate}
\item[{\bf c1)}] Existence and absolute integrability w.r.t. $\mathbf{z}$ of the first two gradients of the LLF $\ln p(\mathbf{z}|\mathbf{x})$, w.r.t. $\mathbf{x}$; and
\item[{\bf c2)}] The support of the LF $p(\mathbf{z}|\mathbf{x})$ be independent of $\mathbf{x}$.
\end{enumerate}

\section{The CRLLB}\label{sec:CRLLB}
It is widely accepted that {\bf c2)} is a necessary condition for the CRLB to hold. In many cases however, the support of the LF depends on ${\bf x}$. For example, this arises if we have the following observation model with additive noise 
\begin{equation}\label{eq:model}
\mathbf{z} = \mathbf{f}(\mathbf{x}) +\mathbf{w}
\end{equation}
where the noise has finite support (e.g., $\mathbf{w}$ has a truncated Gaussian pdf)\footnote{The result to be presented does not necessarily require the model in \eqref{eq:model}; only the LF $p({\bf z}|{\bf x})$ is needed to relate $\mathbf{x}$ to $\mathbf{z}$. However, \eqref{eq:model} is useful to illustrate the dependency of the LF's support on $\mathbf{x}$.}. How can we find a valid CRLB here? To answer this question, we will need to revisit Sec.~II to scrutinize the derivations of CRLB. 


Accounting for parameter-dependent support ${\cal S}({\bf x})$, one can rewrite \eqref{eq:unbiased_CRLB} as
\begin{align}
&\mathbb{E}[(\hat{\mathbf{x}}(\mathbf{z}) - \mathbf{x})|\mathbf{x}] = \int_{\mathcal{S}(\mathbf{x})}p(\mathbf{z}|\mathbf{x})(\hat{\mathbf{x}}(\mathbf{z}) - \mathbf{x})^{\top} d\mathbf{z}= \mathbf{0} \label{eq:unbiased_CRLLB}
\end{align}
where $d\mathbf{z}$ is the infinitesimal hypervolume for the vector $\mathbf{z}$, and is written using the exterior product notation as $d\mathbf{z} = d{z}_{ 1}\wedge d{z}_{2}\wedge...\wedge d {z}_{n_{\mathbf{z}}}$. 

Similarly, we will take the gradient of~\eqref{eq:unbiased_CRLLB} w.r.t. ${\bf x}$. For notational brevity, we will first define the row vector
\begin{align}
\boldsymbol{\Phi}^{\top} \triangleq p(\mathbf{z}|\mathbf{x})(\hat{\mathbf{x}}(\mathbf{z}) - \mathbf{x})^\top
\end{align}
Applying the general Leibniz rule, the resulting gradient, in view of the assumed unbiasedness, is given by 
\begin{align}
&\mathbf{D}(\mathbf{x}) = \nabla_\mathbf{x}\left( \int_{\mathcal{S}({\bf x})} \boldsymbol{\Phi}^{\top} d\mathbf{z}\right) =  \mathbf{D}_{\rm L}(\mathbf{x}) + \mathbf{D}_{\rm I}(\mathbf{x}) = {\bf 0}_{n_{\bf x}}\label{eq:D_1} 
\end{align}
which consists of the contributions from both the integrand 
\begin{align}
& \mathbf{D}_{\rm I}(\mathbf{x}) \triangleq  \int_{\mathcal{S}(\mathbf{x})}\!\!\!\!\nabla_\mathbf{x} \boldsymbol{\Phi}^{\top} d\mathbf{z} = \int_{\mathcal{S}(\mathbf{x})}\!\!\!\!\nabla_\mathbf{x} [p(\mathbf{z}|\mathbf{x})(\hat{\mathbf{x}}(\mathbf{z}) - \mathbf{x})^\top] d\mathbf{z}\label{eq:D_I}
\end{align}
and the support ${\cal S}({\bf x})$ 
\begin{align}\label{eq:D_L_new}
\mathbf{D}_{\rm L}(\mathbf{x})&\triangleq \int_{\delta\mathcal{S}(\mathbf{x})}\left( \nabla_\mathbf{x}\mathbf{z}^{\top} \nabla_{\mathbf{z}}\right) \lrcorner  \left(\boldsymbol{\Phi}^{\top} d\mathbf{z}\right) \;.
\end{align}
In~\eqref{eq:D_L_new}, $\delta\mathcal{S}({\bf x})$ denotes the boundary of $\mathcal{S}({\bf x})$, and the symbol $\lrcorner$ is the interior product operator which, as defined in \eqref{eq:interior product} in Appendix \ref{sec:Leibniz}, performs the degree-1 antiderivation on the exterior product of differential forms in a vector field. Generally speaking, the interior product reduces the dimensionality of the integration by $1$. For example, if $\mathcal{S}({\bf x})$ represents a plane, its boundary $\delta\mathcal{S}({\bf x})$ is a circle (closed line). Therefore, applying the interior product reduces a 2-dimensional integral to 1-dimensional integral. In practice, it is more convenient to use the following formula 
\begin{equation}
\mathbf{D}_{\rm L}(\mathbf{x}) = \int_{\delta\mathcal{S}(\mathbf{x})}\nabla_\mathbf{x}\mathbf{z}^{\top}\, \mathbf{n} p(\mathbf{z}|\mathbf{x})(\hat{\mathbf{x}}(\mathbf{z}) - \mathbf{x})^\top dc \label{eq:D_L_2}
\end{equation}
where $\mathbf{n}$ is the unit vector outward orthogonal to the boundary of the integration volume and $dc$ is the infinitesimal change along the boundary (which is a contour).

Next, we will further decompose  $ \mathbf{D}_{\rm I}(\mathbf{x})$ \eqref{eq:D_I} as 
\begin{align}
\mathbf{D}_{\rm I}(\mathbf{x}) &=  \int_{\mathcal{S}(\mathbf{x})}\!\!\!\!\! [\nabla_\mathbf{x} p(\mathbf{z}|\mathbf{x})](\hat{\mathbf{x}}(\mathbf{z}) - \mathbf{x})^\top d\mathbf{z} - \mathbf{I}_{n_{\mathbf{x}}} \int_{\mathcal{S}(\mathbf{x})}\!\!\!\! p(\mathbf{z}|\mathbf{x})d\mathbf{z} \nonumber\\
& = \mathbf{L}(\mathbf{x}) - \mathbf{I}_{n_{\mathbf{x}}}\label{eq:D_I_c}
\end{align}
where 
\begin{align}\label{eq:Lx}
\mathbf{L}(\mathbf{x}) \triangleq & \int_{\mathcal{S}(\mathbf{x})} \left[\nabla_\mathbf{x} p(\mathbf{z}|\mathbf{x})\right](\hat{\mathbf{x}}(\mathbf{z}) - \mathbf{x})^\top d\mathbf{z}
\end{align}
According to \eqref{eq:D_1} and \eqref{eq:D_I_c}, the following equality holds
\begin{align}
\mathbf{L}(\mathbf{x}) =  \mathbf{I}_{n_{\mathbf{x}}} - \mathbf{D}_{\rm L}(\mathbf{x}) \label{eq:Lx_new}
 \end{align}

Subsequently, we will rely on $\boldsymbol{\epsilon}(\mathbf{z},\mathbf{x})$~\eqref{eq: epsilon} and $\boldsymbol{\gamma}(\mathbf{z},\mathbf{x})$~\eqref{eq: gamma} to express the covariance matrix ${\bf P}({\bf x})$~\eqref{eq:P} and the FIM ${\bf J}({\bf x})$~\eqref{eq:J}. Further leveraging~\eqref{eq:equality}, we can rewrite $\mathbf{L}(\mathbf{x})$ as
\begin{equation}
\mathbf{L}(\mathbf{x}) = \mathbb{E}[\boldsymbol{\gamma}(\mathbf{z},\mathbf{x})\boldsymbol{\epsilon}(\mathbf{z},\mathbf{x})^{\top}|\mathbf{x}]
\end{equation}

Following the proof in \cite{ljung1998system}, we have 
\begin{align}
\mathbb{E}\left[\begin{pmatrix}
\boldsymbol{\epsilon}\\
\boldsymbol{\gamma}
\end{pmatrix} \begin{pmatrix}\boldsymbol{\epsilon}^\top&\boldsymbol{\gamma}^\top \end{pmatrix}\middle\vert\mathbf{x}\right] 
&= 
\begin{bmatrix}
\mathbf{P}(\mathbf{x})&\mathbf{L}(\mathbf{x})^\top\\
\mathbf{L}(\mathbf{x})&\mathbf{J}(\mathbf{x})\end{bmatrix} \geq \mathbf{0}\label{PJL_b}
\end{align}
Therefore, we have the following inequality concerning the covariance matrix $\mathbf{P}(\mathbf{x})$, which is the {\em multidimensional CRLLB}
\begin{subequations}\label{CRLLB}
\begin{align}
\mathbf{P}(\mathbf{x}) &\geq  \mathbf{L}(\mathbf{x})^\top\mathbf{J}(\mathbf{x})^{-1}\mathbf{L}(\mathbf{x})\label{CRLLB_a}\\
                 & =  \left[\mathbf{I}_{n_\mathbf{x}} -  \mathbf{D}_{\rm L}(\mathbf{x})\right]^\top \mathbf{J}(\mathbf{x})^{-1} \left[\mathbf{I}_{n_\mathbf{x}} -  \mathbf{D}_{\rm L}(\mathbf{x})\right]\label{CRLLB_b}
\end{align}
\end{subequations}
where the fact that one matrix is greater than or equal to another matrix means that their difference is positive semi-definite. Eq.\ \eqref{CRLLB_a} follows from Eq.\ \eqref{PJL_b}: when the latter is positive definite, the top left block of the inverse of the partitioned matrix in \eqref{PJL_b} is also positive definite (see, e.g., \cite{BarShalomLK:2001}, Eq.\ (1.3.3-7)) as well; and equality in \eqref{CRLLB_a} --- i.e., statistical efficiency --- holds when \eqref{eq: epsilon} and \eqref{eq: gamma} satisfy the following \textit{generalized collinearity condition}\footnote{This term is used as the multidimensional counterpart of the collinearity condition (see, e.g. \cite{BarShalomLK:2001} Eq. (2.7.3-9) or \cite{VanTrees:2013} Eq. (4-104))  from the scalar CRLB.}
\begin{subequations}\label{collinear}
\begin{align}
\boldsymbol{\epsilon}(\mathbf{z},\mathbf{x}) & = {\bf C}(\mathbf{x}) \boldsymbol{\gamma}(\mathbf{z},\mathbf{x}) \\
\ \hat{\mathbf{x}}(\mathbf{z}) - \mathbf{x} & = {\bf C}(\mathbf{x})\nabla_\mathbf{x} \ln \,p(\mathbf{z}|\mathbf{x})\qquad  a.s.\forall \mathbf{z} \in \mathcal{S}
\end{align}
\end{subequations}
where ``$a.s. \forall \mathbf{z}$" means ``almost surely for all $\mathbf{z}$" according to the probability measure $p(\mathbf{z}|\mathbf{x})$, and ${\bf C}(\mathbf{x})$ is a matrix which depends on $\mathbf{x}$ but not on $\mathbf{z}$, and is given by\footnote{For LFs with parameter-independent support, ${\bf C}(\mathbf{x}) = \mathbf{J}(\mathbf{x})^{-1}$.}
\begin{align}
{\bf C}(\mathbf{x}) = \mathbf{L}(\mathbf{x})^\top \mathbf{J}(\mathbf{x})^{-1}
\end{align}
The above indicates that a necessary condition for the existence of an efficient estimator is that the LF belongs to the exponential family \cite{VanTrees:2013}.

The following remarks are in order.

\noindent{\bf Remark 1.} Note that when the Leibniz term $ \mathbf{D}_{\rm L}(\mathbf{x}) = \mathbf{0}_{n_{\bf x}}$, the CRLLB is the same as the CRLB. Therefore, the CRLLB is a generalization of the CRLB. A sufficient condition for $ \mathbf{D}_{\rm L}(\mathbf{x}) = \mathbf{0}_{n_{\bf x}}$ to hold is that $p({\bf z}|{\bf x}) = 0 (\forall {\bf z} \in \delta {\cal S}({\bf x}))$ (cf.~\eqref{eq:D_L_2}). Thus,
 a sufficient condition for the equivalence of CRLLB and CRLB is that the LF is $0$ at the boundary of its support~\cite{YBS_2014_July}.

\noindent{\bf Remark 2.} The derivations of CRLB in Sec.~II rely on the Cauchy-Schwartz inequality, while the counterparts in the CRLLB employ the positive definiteness of the inverse of a positive definite matrix. The latter is more {\it general}, as it can be used to prove the CRLB as well.

\noindent {\bf Remark 3.} The Leibniz rule for scalar parameters can be regarded as a special case of the multidimensional case. With ${\cal S}({ x}) = [l(x), u(x)]$, the boundary is then $\delta {\cal S}({ x}) = \{l(x), u(x)\}$. When $z = l(x)$, $\nabla_{x} z = l{'} (x)$ and $v_n = -1$. When $z = u(x)$, $\nabla_{x} z = u{'} (x)$ and $v_n = 1$. Then, the Leibniz term is given by
\begin{align}
    L(x) &= u{'} (x) p(z|x)(\hat{x}(z) - x)|_{z = u(x)} \nonumber\\
   & \quad - l{'} (x) p(z|x)(\hat{x}(z) - x)|_{z = l(x)}
\end{align}
which is the same as derived in~\cite{YBS_2014_March} for scalar parameters.

\noindent {\bf Remark 4.} Given the parameter-dependent support, the alternative form of the FIM based on the Hessian of the LLF~\eqref{eq:FIM_hessian} does not hold.

\noindent {\bf Remark 5.} For $M$ i.i.d. measurements, the CRLLB is equal to $1/M$ of their single measurement counterpart.

\section{Measurement model with identity observation matrix}\label{sec:examples}
In this section, we will provide examples of LFs with one vector-valued measurement to demonstrate the usefulness of CRLLB. 

To proceed, consider the measurement model~\eqref{eq:model} with ${\bf f}({\bf x}) = {\bf x}$ ($n_{\bf z} = n_{\bf x} = n$), namely,
\begin{equation}
\mathbf{z} = \mathbf{x} + \mathbf{w} \label{eq:model_X_w}
\end{equation}
where the noise pdf $p_w ({\bf w})$ is a function of 
\begin{align}
  \| {\bf w}\| =  \sqrt{w_1^2 + w_2^2 + \cdots + w_n^2}
\end{align}
with support in an $(n-1)-$sphere\footnote{For any natural number $n$, an $(n-1)$-sphere of radius $a$ is defined as the set of points in $n$-dimensional Euclidean space that are distant from a fixed point $\mathbf{c}$ by $a$ \cite{n_spehre}. In particular: a 0-sphere is a pair of points $\{c-a, c+a \}$; a $1$-sphere is a circle of radius $r$ centered at $\mathbf{c}$; a $2$-sphere is the ordinary sphere in $3$-dimensional Euclidean space, and is the boundary of an ordinary ball.}, namely, ${\cal S}({\bf w}) = \{{\bf w}: \|{\bf w}\| \leq a \}$. Notably, $p({\bf w})$ is assumed to be unimodal with mode at $\|{\bf w}\| = 0$ (see examples later in this section). 

The LF of ${\bf x}$ based on ${\bf z}$ is $p({\bf z}|{\bf x}) = p_w ({\bf z} - {\bf x})$ with support
\begin{align}
    {\cal S}({\bf x}) = \{{\bf x}: \|{\bf x} - {\bf z}\| \leq a \}
\end{align}

Let 
\begin{equation}
r = ||\mathbf{z}-\mathbf{x}||
\end{equation}
and expressing the above in spherical coordinates yields
\begin{equation}
\mathbf{z}-\mathbf{x} = r {\bf n}_n \label{eq:n_sphere}
\end{equation}
where $r \in \left[0,\  a\right]$, and
\begin{align} 
 {\bf n}_n =   \begin{bmatrix}
 \cos\theta_1\\
 \sin\theta_1 \cos\theta_2\\
 \sin\theta_1 \sin\theta_2 \cos\theta_3\\
       \vdots\\
  \sin\theta_1\cdots \sin\theta_{n-2} \cos\theta_{n-1}\\ 
 \sin\theta_1\cdots \sin\theta_{n-2} \sin\theta_{n-1} \label{eq:v_n}
\end{bmatrix}
 \end{align}
with $\theta_1, \theta_2, \cdots, \theta_{n-2} \in \left[0,\ \pi\right]$, $\theta_{n-1} \in \left[0, \ 2\pi\right]$.

The boundary of the support is then given by
\begin{equation}
\delta\mathcal{S}({\bf x}): ||\mathbf{z}-\mathbf{x}|| = a
\end{equation}
Along the boundary, $\mathbf{n}_n$ is the unit vector outwards and normal to the boundary
\begin{align}
&\nabla_\mathbf{x}\mathbf{z}^{\top} = {\bf I}_n
\qquad dc = a d\theta
\end{align}
For this family of LFs, the MLE is given by
\begin{align}
  \hat{\bf x} ({\bf z}) = {\bf z}  \label{eq:MLE}
\end{align}
which is easily seen to be unbiased.

Next, we will give several examples of $p_w$ (LF) for model~\eqref{eq:model_X_w} and the corresponding CRLLB.

\subsection{Raised Fractional Cosine (RFC) within a circle (1-sphere) }
\begin{figure}[!t]
    \centering
    \includegraphics[width=0.5\textwidth]{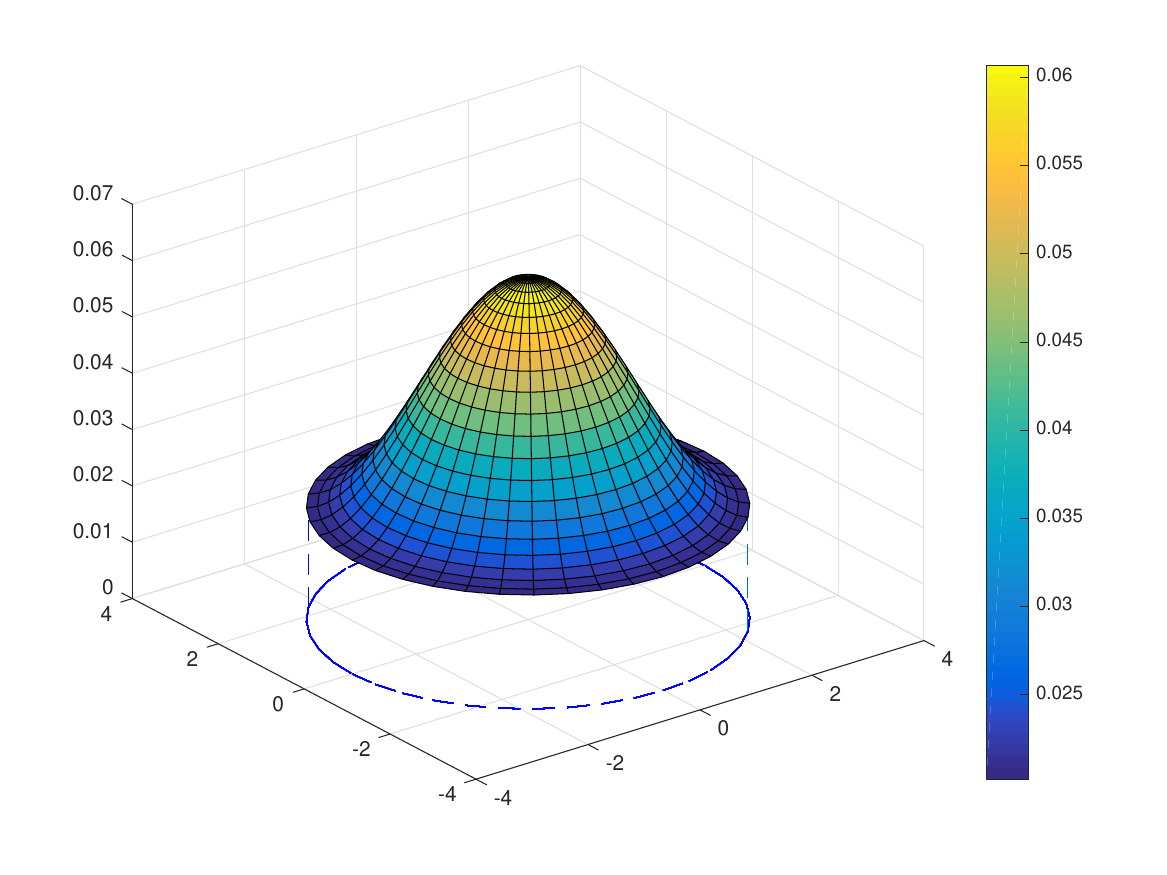}
   \caption{The pdf of RFC measurement noise for $\beta = 0.5$ and $a = \pi$.}
    \label{fig:RFC_beta_05}
\end{figure}
Consider the RFC noise pdf with $n=2$
\begin{equation}\label{eq:pw}
p(\mathbf{w})= \left\{\begin{IEEEeqnarraybox}[\relax][c]{l's}
k_{\rm RFC} \left(1+\beta \cos \left(\frac{\pi ||\mathbf{w}||}{a}\right)\right)&for $||\mathbf{w}||\leq a $\\
0&elsewhere.
\end{IEEEeqnarraybox}\right.
\end{equation}
where $0 \leq \beta \leq 1$ and  $k_{\rm RFC}$ is the normalization constant
\begin{equation}
k_{\rm RFC} = \left(\pi a^2 -4\beta a^2/\pi\right)^{-1}
\end{equation}
Note that the support of $p(\mathbf{w})$ is the region inside a circle of radius $a$ and the boundary is the circle itself. Fig.~\ref{fig:RFC_beta_05} shows an example of the pdf of the measurement noise~\eqref{eq:pw} for $\beta = 0.5$ and $a = \pi$. 


For simplicity, we will set $a = \pi$ in the following. The LF is then given by 
\begin{equation}\label{eq:LF}
p(\mathbf{z}|\mathbf{x})= \left\{\begin{IEEEeqnarraybox}[\relax][c]{l's}
\frac{1+\beta \cos \left(||\mathbf{z}-\mathbf{x}||\right)}{\pi^3 -4\beta \pi}&if $||\mathbf{z}-\mathbf{x}||\leq \pi $\\
0&elsewhere
\end{IEEEeqnarraybox}\right.
\end{equation}
Since the LLF is a function of $r = \|\mathbf{z}-\mathbf{x}\|$, the chain rule yields its gradient w.r.t. $\mathbf{x}$ given by 
\begin{equation}\label{eq:LLF}
\nabla_\mathbf{x} \ln p(\mathbf{z}|\mathbf{x})= \nabla_\mathbf{x} r \frac{d \ln p(\mathbf{z}|\mathbf{x})}{dr}
\end{equation}
where 
\begin{align}
\frac{d \ln p(\mathbf{z}|\mathbf{x})}{dr} & = \frac{-\beta \sin (r)}{1+\beta \cos (r)} \\
\nabla_\mathbf{x} r &= -{\bf n}_2
\end{align}
where ${\bf n}_2$ is as defined in~\eqref{eq:v_n} by setting $n=2$.
The FIM is then 
\begin{align}
\mathbf{J}(\mathbf{x}) &= \int_{\mathcal{S}} \left(\frac{d \rm{\ln}\,p(\mathbf{z}|\mathbf{x})}{d r}\right)^2p(\mathbf{z}|\mathbf{x}) \nabla_\mathbf{x} r \left(\nabla_\mathbf{x} r\right)^{\top} d\mathbf{z}\nonumber\\
&=\frac{1}{\pi^3 -4\beta \pi}\int_{\mathcal{S}} \frac{\beta^2 \sin^2(r)}{1+\beta \cos (r)}\nonumber\\
&\quad \cdot \begin{bmatrix} \cos^2\theta_1 \quad \cos\theta_1 \sin\theta_1\\ \cos\theta_1 \sin\theta_1\quad \sin^2\theta_1 \end{bmatrix}rdr d\theta\nonumber\\
&= \frac{\beta^2}{\pi^2 -4\beta}\int_{r = 0}^{\pi} \frac{ r\sin^2(r)}{1+\beta \cos (r)}dr \mathbf{I}_2\label{eq:J_RFS_b}
\end{align}
which has no analytic expression.

Next, we consider the Leibniz term, which results from the variation of the boundary of the integration area, which is
\begin{equation}
\delta\mathcal{S}({\bf x}): ||\mathbf{z}-\mathbf{x}|| = r = \pi
\end{equation}
 On the boundary, 
 \begin{equation}
 p(\mathbf{z}|\mathbf{x}) = \frac{1-\beta}{\pi^3 -4\beta \pi}
 \end{equation}
 Next, we will calculate the Leibniz term using~\eqref{eq:D_L_2}.  At angle $\theta$ (w.r.t. the x axis) in the circle, we have
\begin{align}
&\mathbf{v}_2 = \begin{bmatrix} \cos\theta\\ \sin\theta \end{bmatrix}
&\nabla_\mathbf{x}\mathbf{z}^{\top} = \begin{bmatrix} 1\quad 0\\ 0\quad 1 \end{bmatrix}
\qquad dc = \pi d\theta
\end{align}
 The Leibniz term in \eqref{eq:D_L_2} becomes
\begin{align}
\mathbf{D}_{\rm L}(\mathbf{x})&= \frac{1-\beta}{\pi^3 -4\beta \pi}\int_{\theta=0}^{2\pi} \begin{bmatrix} 1\quad 0\\ 0\quad 1 \end{bmatrix} \begin{bmatrix} \cos\theta\\ \sin\theta \end{bmatrix} \nonumber\\
&\qquad \cdot \begin{bmatrix} \pi \cos\theta& \pi \sin\theta \end{bmatrix}\pi d\theta \nonumber\\
&=\frac{\pi^2(1-\beta)}{\pi^2 -4\beta}\mathbf{I}_2 \label{eq:D_L_1}
\end{align}
Consequently, $\mathbf{L}(\mathbf{x})$ in \eqref{eq:Lx_new} becomes
\begin{align}
\mathbf{L}(\mathbf{x}) &= \mathbf{I}_2 - \frac{\pi^2(1-\beta)}{\pi^2 -4\beta}\mathbf{I}_2 =\frac{(\pi^2-4)\beta}{\pi^2 -4\beta}\mathbf{I}_2 \label{eq:L_1}
\end{align}
The CRLLB is then 
\begin{align}
\mathbf{P}(\mathbf{x})&\geq\mathbf{L}(\mathbf{x})^\top \mathbf{J}(\mathbf{x})^{-1}\mathbf{L}(\mathbf{x})\nonumber\\
&=\frac{(\pi^2-4)^2}{(\pi^2-4\beta)\int_{r = 0}^{\pi} \frac{ r\sin^2(r)}{1+\beta \cos (r)}dr}\mathbf{I}_2\label{eq:CRLLB_1_b}
\end{align}

\begin{figure}[!t]
    \centering
    \includegraphics[width=0.48\textwidth]{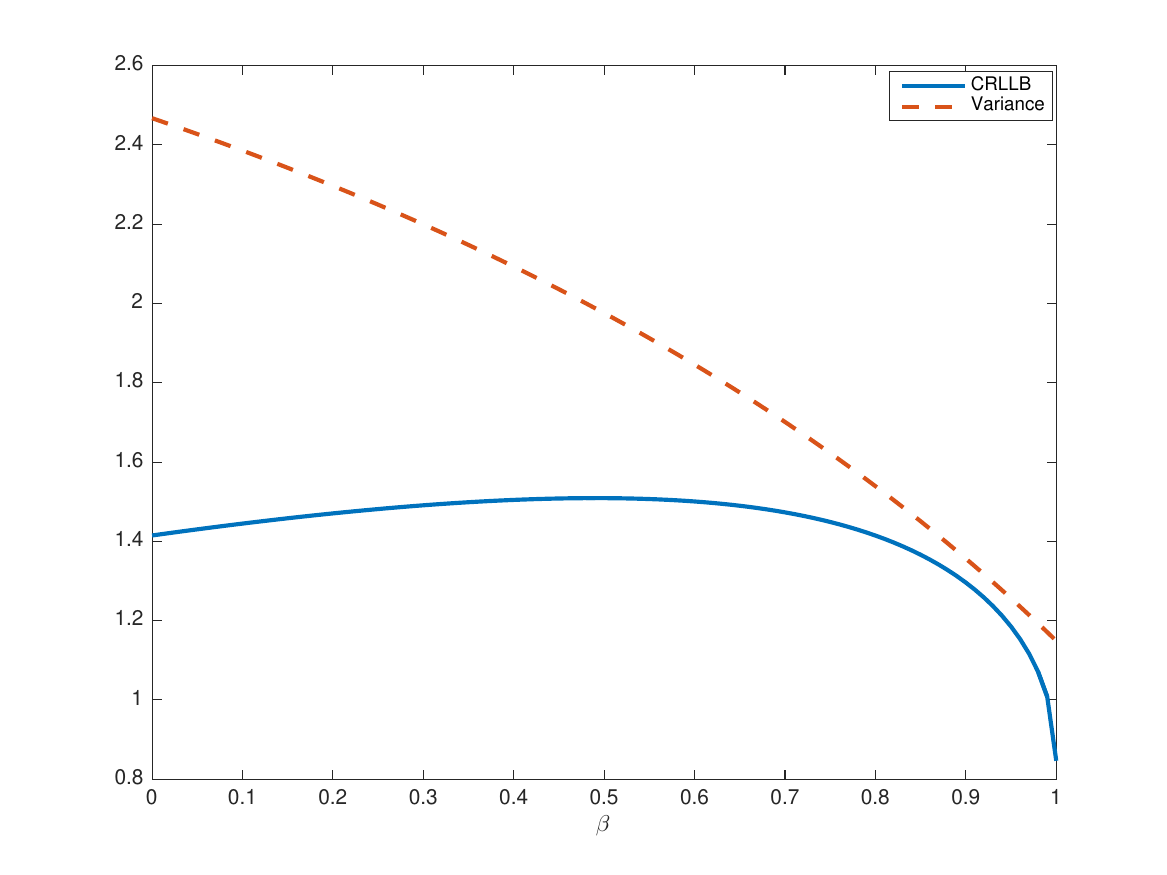}
    \caption{The diagonal terms of the CRLLB vs the counterparts of the covariance of the MLE for a RFC LF with $\beta\in [0,1]$.}
    \label{fig:RC_CRLLB_Var}
\end{figure}

By comparison, the covariance of \eqref{eq:MLE} is
\begin{align}
{\rm cov}\left[\hat{\mathbf{x}}(\mathbf{z})\right]& = \int_{\mathcal{S}({\bf x})}p(\mathbf{z}|\mathbf{x})(\mathbf{z} - \mathbf{x})(\mathbf{z} - \mathbf{x})^\top d\mathbf{z}\nonumber\\
\begin{split}
&=\frac{1}{\pi^3-4\beta\pi}\int_{\theta =0}^{2\pi}\int_{r = 0}^{\pi} (1+\beta \cos r) \nonumber\\
&\quad\cdot \begin{bmatrix} r^2 \cos^2\theta \quad r^2 \cos\theta \sin\theta\\ r^2 \cos\theta \sin\theta\quad r^2 \sin^2\theta \end{bmatrix}rdr d\theta
\end{split}\nonumber\\
&=\frac{1}{\pi^2-4\beta}\left[\frac{\pi^4}{4}+\beta(12-3\pi^2)\right]\mathbf{I}_2\label{eq:var_RFC}
\end{align}
Fig.~\ref{fig:RC_CRLLB_Var} shows the plots of the diagonal terms of \eqref{eq:CRLLB_1_b} and \eqref{eq:var_RFC} for different values of $\beta$. Clearly, the LF \eqref{eq:LF} does not belong to the exponential family, so there is no efficient estimator.

\subsection{Truncated Laplace distribution inside a circle}
Here, we consider the $2$-D measurement noise has a truncated Laplace distribution inside a circle, i.e.,
\begin{equation}\label{eq:pw_laplace}
p({\bf w})= \left\{\begin{IEEEeqnarraybox}[\relax][c]{l's}
k_{\rm TL} \exp (-\alpha ||{\bf w}||)&for $||{\bf w}||\leq a $\\
0&elsewhere
\end{IEEEeqnarraybox}\right.
\end{equation}
where $k_{\rm TL} $ is the normalization constant given by
\begin{equation}
k_{\rm TL}  = \frac{\alpha^2}{2\pi (1-e^{-a\alpha}-a\alpha e^{-a\alpha})}
\end{equation}
For a single measurement, the LF is
\begin{equation}\label{LF_Trunc_Laplace}
p({\bf z}|{\bf x})=\left\{\begin{IEEEeqnarraybox}[\relax][c]{l's}
k_{\rm TL}  \exp \left(-\alpha ||{\bf z}-{\bf x}|| \right) & for $||{\bf z}-{\bf x}||\leq a $\\
0&elsewhere
\end{IEEEeqnarraybox}\right.
\end{equation}
The gradient of the LLF $\rm{\ln}\,p({\bf z}|{\bf x})$ as in \eqref{eq:LLF} is given by
\begin{align}
\nabla_{\bf x} \ln p({\bf z}|{\bf x})&=  \nabla_{\bf x} r \frac{d \ln p({\bf z}|{\bf x})}{d r} = \begin{bmatrix} \alpha\cos\theta \\ \alpha\sin\theta \end{bmatrix} \label{eq:LLF_laplace}
\end{align}
Note the above gradient does not exist at the origin. But it is a.s. differentiable because the point where it is not differentiable has probability zero and thus it does not contribute to the (Lebesgue) integral that yields the FIM, which is given by
\begin{align}
\mathbf{J}({\bf x}) &= k_{\rm TL}  \int_{\theta =0}^{2\pi}\int_{r = 0}^{a}  \exp(-\alpha r)\alpha^2\nonumber\\
&\quad \cdot \begin{bmatrix} \cos^2\theta \quad \cos\theta \sin\theta\\ \cos\theta \sin\theta\quad \sin^2\theta \end{bmatrix}rdr d\theta\nonumber\\
&= \frac{\alpha^2}{2} \mathbf{I}_2 \label{eq:J_Laplace}
\end{align}
The Leibniz term using \eqref{eq:Leibniz_2D_3D} is
\begin{align}\label{eq:D_L_laplace}
\mathbf{D}_{\rm L}(\mathbf{x})& = k_{\rm TL}\!\!  \int_{\theta=0}^{2\pi} \begin{bmatrix} 1\quad 0\\ 0\quad 1 \end{bmatrix} \begin{bmatrix} \cos\theta\\ \sin\theta \end{bmatrix} \begin{bmatrix}  a\cos\theta \!& \! a\sin\theta \end{bmatrix}e^{-a\alpha} a d\theta\nonumber\\
& =  k_{\rm TL} \pi a^2 e^{-a\alpha}\mathbf{I}_2
\end{align}
Then
\begin{align}\label{eq:L_laplace}
\mathbf{L}(\mathbf{x})& = \mathbf{I}_2 - \mathbf{D}_{\rm L}(\mathbf{x}) \nonumber\\
& = \frac{2-2e^{-a\alpha}-2a\alpha e^{-a\alpha} - a^2\alpha^2 e^{-a\alpha}}{2-2e^{-a\alpha}-2a\alpha e^{-a\alpha}}\mathbf{I}_2
\end{align}
The CRLLB is given by
\begin{align}\label{eq:CRLLB_laplace}
\mathbf{P}({\bf x})&\geq  \mathbf{L}({\bf x})^T\mathbf{J}({\bf x})^{-1}\mathbf{L}({\bf x})\nonumber \\
&= \frac{(2-2e^{-a\alpha}-2a\alpha e^{-a\alpha} - a^2\alpha^2 e^{-a\alpha})^2}{2\alpha^2(1-e^{-a\alpha} -a\alpha e^{-a\alpha})^2}\mathbf{I}_2
\end{align}
The covariance of the MLE is
\begin{align}\label{eq:var_laplace}
&{\rm cov}\left[\hat{{\bf x}}({\bf z})\right]= \int_{\mathcal{S}({\bf x})}p({\bf z}|{\bf x})({\bf z} - {\bf x})({\bf z} - {\bf x})^T d{\bf z} \nonumber\\
&= k_{\rm TL}\int_{\theta =0}^{2\pi}\int_{r = 0}^{a} e^{-\alpha r} \begin{bmatrix} r^2 \cos^2\theta \quad r^2 \cos\theta \sin\theta\\ r^2 \cos\theta \sin\theta\quad r^2 \sin^2\theta \end{bmatrix}rdr d\theta \nonumber\\
&=\frac{6 - 6e^{-a\alpha} - 6a\alpha e^{-a\alpha} - 3a^2\alpha^2 e^{-a\alpha} - a^3\alpha^3 e^{-a\alpha}}{2\alpha^2(1-e^{-a\alpha} -a\alpha e^{-a\alpha})} \mathbf{I}_2
\end{align}
\begin{figure}[!t]
    \centering
    \includegraphics[width=0.48\textwidth]{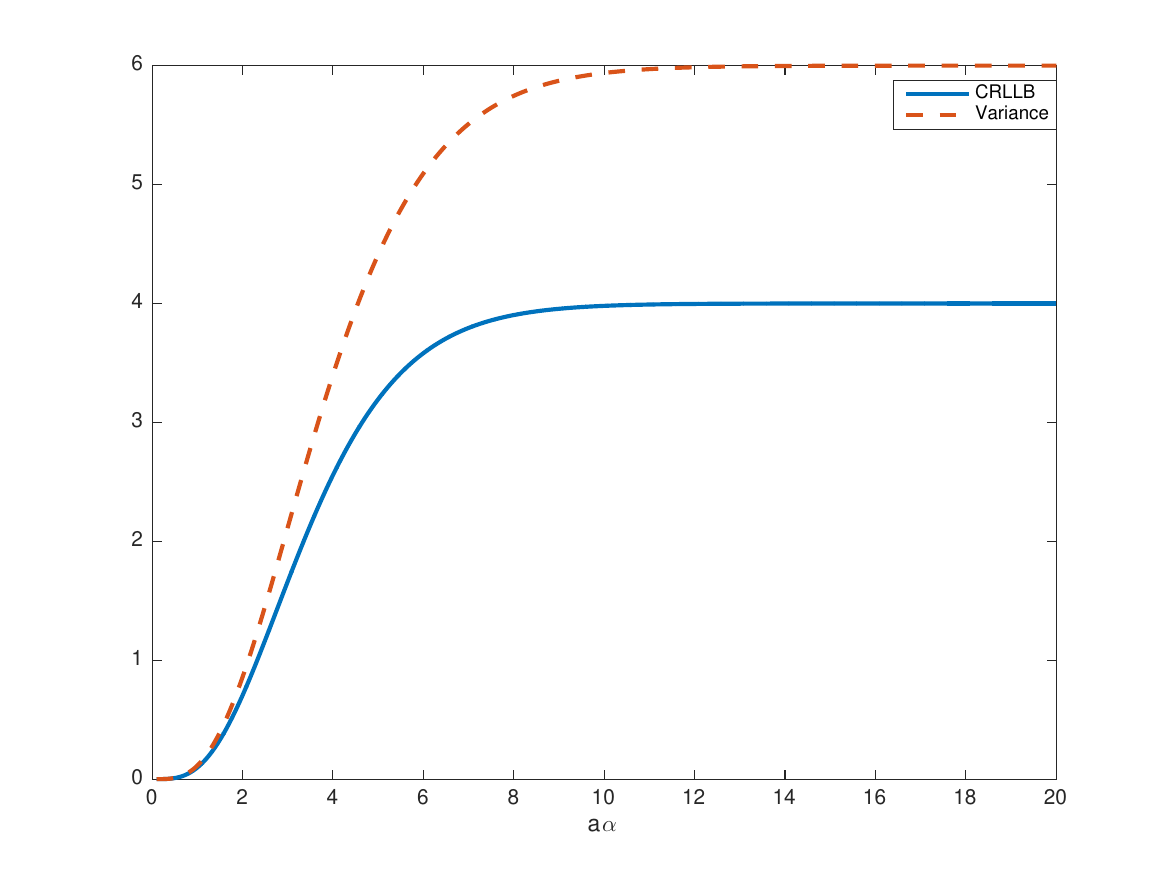}
    \caption{The diagonal terms of the CRLLB vs that of the covariance matrix of the MLE from a truncated Laplace LF for different $a\alpha$.}
    \label{fig:LP_CRLLB_Var}
\end{figure}
Fig. \ref{fig:LP_CRLLB_Var} shows the diagonal terms of the covariance \eqref{eq:var_laplace} and the CRLLB \eqref{eq:CRLLB_laplace} as a function of $a\alpha$ after removing the common factor $(2\alpha^2(1-e^{-a\alpha} -a\alpha e^{-a\alpha}))^{-1}$. Although the LF \eqref{eq:pw_laplace} belongs to the exponential family, the unbiased MLE estimator \eqref{eq:MLE} is not efficient as the \textit{generalized collinearity condition} \eqref{collinear} is not satisfied. This validates that the exponential form of the LF is only a necessary condition for an efficient estimator to hold.

\subsection{Truncated Gaussian (TG) distribution}
Here, the $n$-D measurement noise has a TG distribution in an $(n-1)$-sphere with radius $a$, i.e.,
\begin{equation}\label{eq:pw_TG_nD}
p({\bf w})= \left\{\begin{IEEEeqnarraybox}[\relax][c]{l's}
\frac{1}{k^n_{\rm TG}}\exp(-\frac{||{\bf w}||^2}{2\sigma^2})&for $||{\bf w}||\leq a $\\
0&elsewhere
\end{IEEEeqnarraybox}\right.
\end{equation}
where the normalization constant $k^n_{\rm TG}$ is obtained by integrating \eqref{eq:pw_TG_nD} over the volume of the $n$-sphere as
\begin{align}\label{eq:k_n}
&k_{\rm TG}^n = \int_{\theta_{n-1} = 0}^{2\pi} \int_{\theta_{n-2} = 0}^{\pi} \cdots \int_{\theta_{1} = 0}^{\pi} \int_{r = 0}^{a} \exp(-\frac{r^2}{2})r^{n-1} \nonumber\\ 
&\quad \cdot \sin^{n-2} \theta_1\cdots \sin\theta_{n-2} dr d\theta_1\cdots d\theta_{n-2} d\theta_{n-1} \nonumber\\
& = 2\pi\left(\int_{r = 0}^{a} \exp\left(-\frac{r^2}{2\sigma^2}\right)r^{n-1} dr\right) \prod_{i = 1}^{n-2} \left(\int_{\theta_{i} = 0}^{\pi}\sin^{n-1-i} \theta_i d\theta_i\right) \nonumber
\end{align}

The LF is then given by
\begin{equation}\label{LF_Trunc_Gaussian_nD}
p(\mathbf{z}|\mathbf{x})= \left\{\begin{IEEEeqnarraybox}[\relax][c]{l's}
\frac{1}{k^n_{\rm TG}}\exp (-\frac{||\mathbf{z}-\mathbf{x}||^2}{2\sigma^2})&for $||{\bf z}-{\bf x}||\leq a $\\
0&elsewhere
\end{IEEEeqnarraybox}\right.
\end{equation}
It is easy to verify that the MLE~\eqref{eq:MLE} $\hat{\bf x}({\bf z}) = {\bf z}$ is unbiased. Using the spherical coordinates~\eqref{eq:n_sphere} to express ${\bf z} - {\bf x}$, the 
 covariance of the MLE is then
\begin{align}
&{\rm cov}\left[\hat{\bf x}(\bf z)\right]=  \int_{\mathcal{S}}p({\bf z}|{\bf x})({\bf z} - {\bf x})({\bf z} - {\bf x})^\top d{\bf z} \nonumber\\
& = \left(k^n_{\rm TG}\right)^{-1} \int_{\theta_{n-1} = 0}^{2\pi} \int_{\theta_{n-2} = 0}^{\pi} \!\!\cdots \int_{\theta_{1} = 0}^{\pi} \int_{r = 0}^{a} \exp(-r^2/2) r^{n+1} \nonumber\\
&\cdot \begin{bmatrix} 
\cos^2 \theta_1 & \cdots &  \cos\theta_1\sin\theta_1\cdots \sin\theta_{n-1} \\
\vdots &\ddots &\vdots \\
\cos\theta_1\sin\theta_1\cdots \sin\theta_{n-1}&\cdots  & \sin^2\theta_1\cdots \sin^2\theta_{n-1} \\
  \end{bmatrix}\nonumber\\
  &\cdot\sin^{n-2} \theta_1\cdots \sin\theta_{n-2} dr d\theta_1\cdots d\theta_{n-2} d\theta_{n-1}\nonumber\\
& = \left(k^n_{\rm TG}\right)^{-1} \left(\int_{r = 0}^{a} \exp(-\frac{r^2}{2\sigma^2})r^{n+1} dr\right)  \left(\int_{\theta_{1} = 0}^{\pi}\sin^{n} \theta_1 d\theta_1\right) \nonumber\\
&\quad \cdots\left(\int_{\theta_{n-1} = 0}^{2\pi} \sin^2\theta_{n-1} d\theta_{n-1} \right) \mathbf{I}_n \label{eq:cov_TG_nD_1}
\end{align}
Making use of \eqref{eq:integral_equality_1}--\eqref{eq:integral_equality_3} in Appendix \ref{sec:useful} yields
\begin{align}
&{\rm cov}\left[\hat{\bf x}({\bf z})\right] = \nonumber\\
&\frac{\pi\sigma^2\left(-a^n\exp(-\frac{a^2}{2}) +  n \int_{r = 0}^{a} \exp(-\frac{r^2}{2\sigma^2})r^{n-1} dr\right) }{2\pi \int_{r = 0}^{a} \exp(-\frac{r^2}{2\sigma^2})r^{n-1} dr } \nonumber\\
&\cdot \frac{ \frac{n-1}{n}\left(\int_{\theta_{1} = 0}^{\pi}\sin^{n-2} \theta_1 d\theta_1\right) \cdots\left(\frac{2}{3} \int_{\theta_{n-2} = 0}^{\pi} \sin\theta_{n-2}d\theta_{n-2}\right)}{ \left(\int_{\theta_{1} = 0}^{\pi}\sin^{n-2} \theta_1 d\theta_1\right) \cdots\left( \int_{\theta_{n-2} = 0}^{\pi} \sin\theta_{n-2}d\theta_{n-2}\right)}\mathbf{I}_n \nonumber\\
& = \sigma^2\left(1 - \frac{a^n\exp(-\frac{a^2}{2\sigma^2}) }{n\int_{r = 0}^{a} \exp(-\frac{r^2}{2\sigma^2})r^{n-1} dr}\right) \mathbf{I}_n\label{eq:cov_TG_nD}
\end{align}

\begin{figure}[!t]
    \centering
    \includegraphics[width=0.48\textwidth]{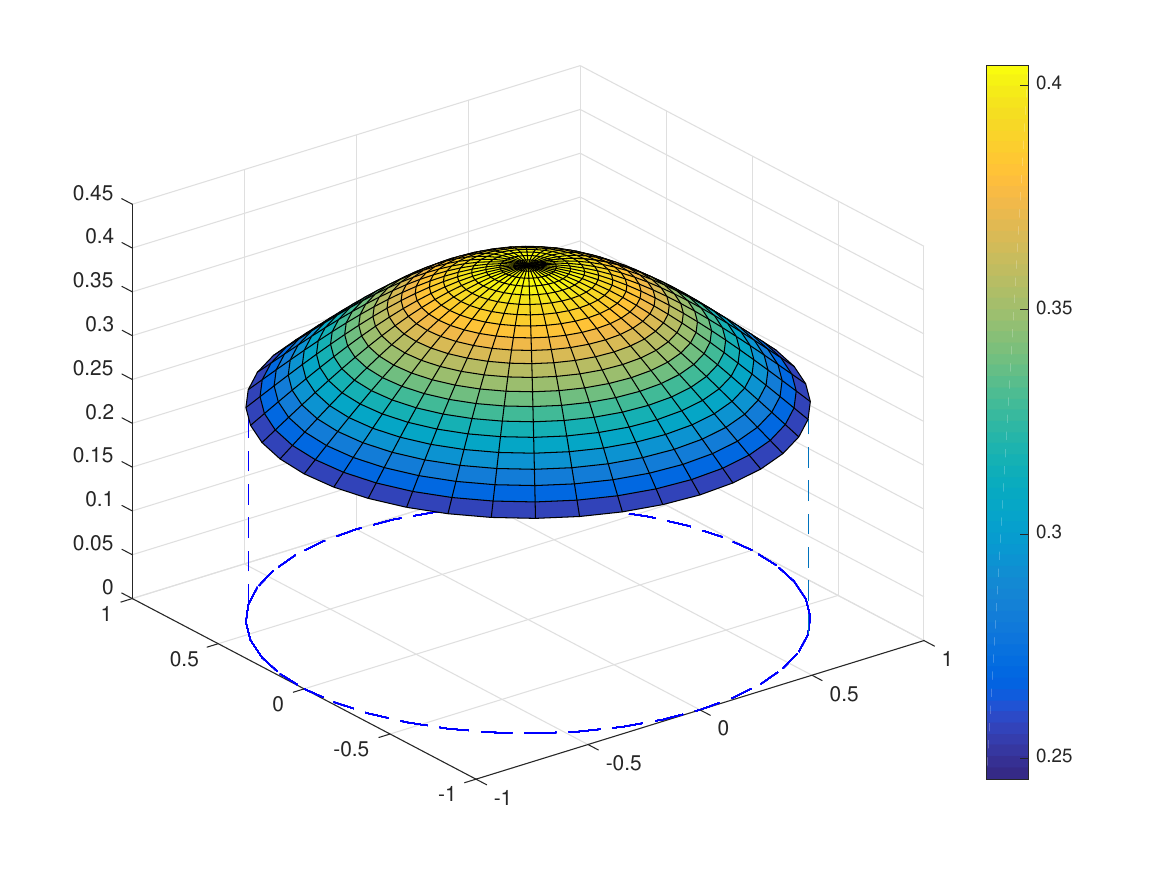}
    \caption{The pdf of TG measurement noise for $n = 2$, $\sigma=1$, and $a = 1$.}
    \label{fig:LP_CRLLB_Var}
\end{figure}

Next, we will evaluate the CRLLB. Firstly, the gradient of the LLF $\ln p({\bf z}|{\bf x})$ according to \eqref{eq:LLF} is given by
\begin{align}
\nabla_{\bf x} \ln p({\bf z}|{\bf x})&= \nabla_{\bf x} r \frac{d \ln p({\bf z}|{\bf x})}{d r} \nonumber \\
&=\begin{bmatrix}
 \cos\theta_1\\
 \sin\theta_1 \cos\theta_2\\
 \sin\theta_1 \sin\theta_2 \cos\theta_3\\
       \vdots\\
  \sin\theta_1\cdots \sin\theta_{n-2} \cos\theta_{n-1}\\ 
 \sin\theta_1\cdots \sin\theta_{n-2} \sin\theta_{n-1}
\end{bmatrix} (-r/\sigma^2) \label{eq:LLF_TG_nD}
\end{align}
The FIM is then 
\begin{align}
\mathbf{J}({\bf x}) & = \int_{\mathcal{S}({\bf x})}p({\bf z}|{\bf x})\nabla_{\bf x} \ln p({\bf z}|{\bf x}) \nabla^{\top}_{\bf x} \ln p({\bf z}|{\bf x}) d{\bf z} \nonumber\\
& = {\rm cov}\left[\hat{{\bf x}}({\bf z})\right]/\sigma^4\label{eq:J_TG_nD}
\end{align}
Evaluating~\eqref{eq:D_L_2} yields the Leibniz term given by
\begin{align}
&\mathbf{D}_{\rm L}({\bf x})  = \left(k^n_{\rm TG}\right)^{-1} \int_{\theta_{n-1} = 0}^{2\pi} \int_{\theta_{n-2} = 0}^{\pi} \cdots \int_{\theta_{1} = 0}^{\pi}
\begin{bmatrix} 
1& \cdots & 0 \\ \vdots& \ddots&\vdots\\ 0 & \cdots& 1
 \end{bmatrix}\nonumber\\
 &\cdot\begin{bmatrix}
 \cos\theta_1\\
       \vdots\\
\sin\theta_1\cdots \sin\theta_{n-1}
\end{bmatrix} 
\begin{bmatrix}
 \cos\theta_1& \cdots &
\sin\theta_1\!\cdots \! \sin\theta_{n-1}
\end{bmatrix}  \nonumber\\
&\cdot e^{-a^2/(2\sigma^2)} a^{n}  \sin^{n-2} \theta_1\cdots \sin\theta_{n-2}\  dr d\theta_1\cdots d\theta_{n-2} d\theta_{n-1}\nonumber\\
& = \frac{a^n\exp(-\frac{a^2}{2\sigma^2}) }{n\int_{r = 0}^{a} \exp(-\frac{r^2}{2\sigma^2})r^{n-1} dr}\mathbf{I}_n \label{eq:D_L_TG_nD}
\end{align}

Then $\mathbf{L}({\bf x})$ is
\begin{align}
\mathbf{L}({\bf x})&= \mathbf{I}_n - \frac{a^n\exp(-\frac{a^2}{2}) }{n\int_{r = 0}^{a} \exp(-\frac{r^2}{2})r^{n-1} dr}\mathbf{I}_n ={\rm cov}\left[\hat{\bf x}({\bf z})\right]/\sigma^2 \label{eq:L_TG_nD}
\end{align}

Finally, the CRLLB is given by
\begin{align}
\mathbf{P}({\bf x})&\geq \mathbf{L}({\bf x})^\top\mathbf{J}({\bf x})^{-1}\mathbf{L}({\bf x})= {\rm cov}\left[\hat{\bf x}({\bf z})\right]\label{eq:CRLLB_TG_nD_b}
\end{align}
which is the same as its covariance \eqref{eq:cov_TG_nD}. This statistical efficiency results from the fact that  the \textit{generalized collinearity condition} \eqref{collinear} is satisfied. Note that \eqref{eq:CRLLB_TG_nD_b} is actually independent of $\mathbf{x}$.

\subsection{Uniform distribution in a circle}
For the TG noise pdf~\eqref{eq:pw_TG_nD}, it represents a uniform distribution within in a circle with $n = 2$ and $\sigma\rightarrow \infty$. When $n = 2$, the covariance matrix~\eqref{eq:cov_TG_nD} of the MLE is given by
\begin{subequations}
\begin{align}
 &{\rm cov}\left[\hat{\bf x}(\bf z)\right]= \sigma^2\left(1- \frac{\frac{a^2}{2\sigma^2}\exp(-\frac{a^2}{2\sigma^2})}{1 - \exp(-\frac{a^2} {2\sigma^2})}\right)  {\bf I}_2 \nonumber\\
 & = \frac{\sigma^2[1- (1 + \frac{a^2}{2\sigma^2})\exp(-\frac{a^2}{2\sigma^2})]}{1 - \exp(-\frac{a^2}{2\sigma^2})}{\bf I}_2 \nonumber\\
 & \approx \frac{\sigma^2[1- (1 + \frac{a^2}{2\sigma^2})\exp(-\frac{a^2}{2\sigma^2})]}{1 - (1-\frac{a^2}{2\sigma^2})} {\bf I}_2\label{eq:approx_1}\\
 & \approx \frac{2\sigma^4}{a^2}\left[1 -\left(1 + \frac{a^2}{2\sigma^2}\right)\left(1 - \frac{a^2}{2\sigma^2} + \frac{a^4}{8\sigma^4}\right) \right]{\bf I}_2 \label{eq:approx_2} \\
  &=  \frac{a^2}{4}{\bf I}_2 \qquad \sigma\rightarrow \infty    \label{eq:approx_3}
\end{align}
\end{subequations}
where $\exp(-\frac{a^2}{2\sigma^2})\approx 1 - \frac{a^2}{2\sigma^2}$ in the denominator yields the approximation in~\eqref{eq:approx_1}; the approximation in~\eqref{eq:approx_2}  comes from $\exp(-\frac{a^2}{2\sigma^2})\approx 1 - \frac{a^2}{2\sigma^2} + \frac{a^4}{8\sigma^4}$ in the numerator; and ignoring the high-order terms as $\sigma\rightarrow \infty$ yields~\eqref{eq:approx_3}. 

As shown in~\eqref{eq:CRLLB_TG_nD_b}, the CRLLB for the uniform LF using the TG approximation is the same as the covariance $\frac{a^2}{4}$. On the other hand, the uniform distribution can also be given by  the RFC pdf
\eqref{eq:pw} with $\beta = 0$ in. The corresponding CRLLB is obtained by evaluating $\beta=0$ in~\eqref{eq:CRLLB_1_b} (for $a = \pi$)
\begin{align}
\mathbf{P}(\mathbf{x})&\geq \mathbf{L}(\mathbf{x})^{\top}\mathbf{J}(\mathbf{x})^{-1}\mathbf{L}(\mathbf{x})\nonumber \\
& =  \frac{(\pi^2-4)^2}{\pi^2\int_{r = 0}^{\pi} r\sin^2(r)dr}\mathbf{I}_2 = 1.4147\mathbf{I}_2 \label{eq:CRLLB_2}
\end{align}
The covariance of the estimator for this case of uniformly distributed noise inside a circle is obtained by setting $\beta$ to 0 in \eqref{eq:var_RFC} or by setting $a = \pi$ in~\eqref{eq:approx_3}
\begin{equation}\label{eq:var_1}
{\rm cov}\left[\hat{\mathbf{x}}(\mathbf{z})\right]= 2.4674\mathbf{I}_2
\end{equation}
It is interesting to see that different approximations yield different CRLLBs. Further, the claim of statistical efficiency varies depending on the form of the approximated LF -- the TG-based approximation yields an efficient estimator, while the RFC-based counterpart does not.

\section{Measurement model with non-identity observation matrix}\label{sec:non-identity}
The examples in the previous section dealt with the case of identity observation matrix. In this section, we provide an example with non-identity observation matrix with truncated Gaussian distribution for the measurement noise, and show that the CRLLB still applies. Consider the following observation model\footnote{This is an illustration of the case where $n_{\mathbf{x}} \neq n_{\mathbf{z}}$, and the approach holds for a general $n_{\mathbf{z}} > n_{\mathbf{x}}$ for a single measurement.}
\begin{align}
\mathbf{z} = \mathbf{H}\mathbf{x} + \mathbf{w}
\end{align}
where $\mathbf{x}$ is the $2$-dimensional vector to be estimated ($n_{\mathbf{x}} = 2$), $\mathbf{z}$ is the $3$-dimensional measurement ($n_{\mathbf{z}} = 3$), $\mathbf{H}$ is the $3\times 2$ observation matrix with rank $2$ and $\mathbf{w}$ is the $3$-dimensional measurement noise with a TG pdf inside a 2-sphere as in~\eqref{eq:pw_TG_nD} with $n=3$ and the normalization constant given by
\begin{equation}
k^3_{\rm TG} = \frac{1}{4\pi \left(\int_{r = 0}^{a} e^{-r^2/2}dr - ae^{-r^2/2} \right)}
\end{equation}

The LF is given by
\begin{align}\label{LF_Trunc_Gaussian_3D}
p(\mathbf{z}|\mathbf{x})= \left\{\begin{IEEEeqnarraybox}[\relax][c]{l's}
k^3_{\rm TG} e^{-\frac{\left(\mathbf{z}-\mathbf{H}\mathbf{x}\right)^{\top}\left(\mathbf{z}-\mathbf{H}\mathbf{x}\right)}{2}} & for $||{\bf z}-{\bf x}||\leq a $\\
0&elsewhere
\end{IEEEeqnarraybox}\right.
\end{align}
The MLE is then given by
\begin{align}\label{eq:MLE_TG_H}
\hat{\mathbf{x}}(\mathbf{z}) = \mathbf{H}^{\dag}\mathbf{z}
\end{align}
where $ \mathbf{H}^{\dag}$ is the pseudo inverse as
\begin{align}
 \mathbf{H}^{\dag} = \left( \mathbf{H}^{\top} \mathbf{H}\right)^{-1}  \mathbf{H}^{\top}
\end{align}
As
\begin{align}
\mathbb{E}\left( \hat{\mathbf{x}}(\mathbf{z}) -  \mathbf{x}\right) = \mathbb{E}\left(  \mathbf{H}^{\dag}\mathbf{w} \right) = \mathbf{0}
\end{align}
the MLE \eqref{eq:MLE_TG_H} is also an unbiased estimator.

Let 
\begin{equation}
r = ||\mathbf{z}-\mathbf{H}\mathbf{x}||
\end{equation}
expressing in spherical coordinates yields
\begin{equation}
\mathbf{z}-\mathbf{H}\mathbf{x} = r{\bf n}_3 = 
\begin{bmatrix}
r\cos\theta_1\\
r\sin\theta_1\cos\theta_2\\
r\sin\theta_1\sin\theta_2
\end{bmatrix}
\end{equation}
where $\theta_1 \in \left[0,\ \pi\right]$ is the elevation angle, and $\theta_{2} \in \left[0, \ 2\pi\right]$ is the azimuth angle.

The covariance of the estimator \eqref{eq:MLE_TG_H} is then
\begin{align}
&{\rm cov}\left[\hat{\bf x}({\bf z})\right]= \int_{\mathcal{S}({\bf x})}p(\mathbf{z}|\mathbf{x})(\mathbf{H}^{\dag} \mathbf{z} - \mathbf{x})(\mathbf{H}^{\dag} \mathbf{z} - \mathbf{x})^{\top} d\mathbf{z} \nonumber\\
& = \mathbf{H}^{\dag} \int_{\mathcal{S}({\bf x})}p(\mathbf{z}|\mathbf{x})(\mathbf{z} - \mathbf{H}\mathbf{x})(\mathbf{z} - \mathbf{H}\mathbf{x})^\top d\mathbf{z} \left(\mathbf{H}^{\dag}\right)^{\top} \nonumber\\
& =\mathbf{H}^{\dag} k^3_{\rm TG} \int_{\theta_1= 0}^{\pi}\int_{\theta_2 =0}^{2\pi}\int_{r = 0}^{a} \exp(-r^2/2) \nonumber\\
  &\cdot \begin{bmatrix} 
\cos^2\theta_1\! &\! \cos\theta_1\!\sin\theta_1\!\cos\theta_2 \!&\!\cos\theta_1\!\sin\theta_1\!\sin\theta_2\! \\
\cos\theta_1\!\sin\theta_1\!\cos\theta_2\!  &\! \sin^2\theta_1\!\cos^2\theta_2\! &\! \sin^2\theta_1\! \cos\theta_2 \!\sin\theta_2\!  \\
\cos\theta_1\!\sin\theta_1\!\sin\theta_2 \!& \! \sin^2\theta_1\! \cos\theta_2 \!\sin\theta_2\! &\! \sin^2\theta_1\!\sin^2\theta_2 \!
  \end{bmatrix}\nonumber\\  
  &\cdot r^4 \sin\theta_1 dr d\theta_1 d\theta_2 \left(\mathbf{H}^{\dag}\right)^{\top}  \nonumber\\
& = \left(1 -  \frac{4\pi k^3_{\rm TG}}{3} a^3 e^{-a^2/2} \right) \left( \mathbf{H}^{\top} \mathbf{H}\right)^{-1}\label{eq:cov_TG_3D}
\end{align}
Next, we evaluate the CRLLB. Firstly, the gradient of the LLF $\ln p({\bf z}|{\bf x})$ according to \eqref{eq:LLF} is given by
\begin{align}
\nabla_{\bf x} \ln p({\bf z}|{\bf x})& = \nabla_{\mathbf{x}} r \frac{d \ln p(\mathbf{z}|\mathbf{x})}{d r} \nonumber \\
&= \mathbf{H}^{\top} \begin{bmatrix}
-\cos\theta_1\\
-\sin\theta_1\cos\theta_2\\
-\sin\theta_1\sin\theta_2
\end{bmatrix}(-r) \label{eq:LLF_TG_3D}
\end{align}
The FIM is then 
\begin{align}
&\mathbf{J}({\bf x}) = k^3_{\rm TG} \int_{\theta_1= 0}^{\pi}\int_{\theta_2 =0}^{2\pi}\int_{r = 0}^{a} (-r)^2 \exp(-r^2/2)\mathbf{H}^{\top}\nonumber \\
 &\cdot \begin{bmatrix} 
\cos^2\theta_1\! &\! \cos\theta_1\!\sin\theta_1\!\cos\theta_2 \!&\!\cos\theta_1\!\sin\theta_1\!\sin\theta_2\! \\
\cos\theta_1\!\sin\theta_1\!\cos\theta_2\!  &\! \sin^2\theta_1\!\cos^2\theta_2\! &\! \sin^2\theta_1\! \cos\theta_2 \!\sin\theta_2\!  \\
\cos\theta_1\!\sin\theta_1\!\sin\theta_2 \!& \! \sin^2\theta_1\! \cos\theta_2 \!\sin\theta_2\! &\! \sin^2\theta_1\!\sin^2\theta_2 \!
  \end{bmatrix}\nonumber\\  
  &\cdot \mathbf{H} r^2 \sin\theta_1 dr d\theta_1 d\theta_2 \nonumber \\
& = \left(1 -  \frac{4\pi k^3_{\rm TG}}{3} a^3 e^{-a^2/2} \right)\mathbf{H}^{\top}\mathbf{H} \label{eq:J_TG_3D}
\end{align}
For the Leibniz term, using the simplification given in \eqref{eq:Leibniz_2D_3D}, we have 
\begin{align}
\mathbf{D}_{\rm L}({\bf x})& = k^3_{\rm TG}\int_{\theta_1= 0}^{\pi}\int_{\theta_2 =0}^{2\pi} \mathbf{H}^{\top}
\begin{bmatrix}
 \cos\theta_1\\ \sin\theta_1\cos\theta_2\\ \sin\theta_1\sin\theta_2
\end{bmatrix} \nonumber \\
&\quad \cdot \begin{bmatrix}
a\cos\theta_1&a\sin\theta_1\cos\theta_2& a\sin\theta_1\sin\theta_2
\end{bmatrix}e^{-a^2/2} \nonumber\\
&\quad  \cdot \left(\mathbf{H}^{\dag}\right)^{\top} a^2 \sin\theta_1d\theta_1 d\theta_2 \nonumber\\
&=\frac{4\pi k^3_{\rm TG}}{3} a^3 e^{-a^2/2}\mathbf{I}_3 \label{eq:D_L_TG_3D}
\end{align}
Then $\mathbf{L}({\bf x})$ is
\begin{align}
\mathbf{L}({\bf x})& = \mathbf{I}_3 - \frac{4\pi k^3_{\rm TG}}{3} a^3 e^{-a^2/2}\mathbf{I}_3 =\left(1 -  \frac{4\pi k^3_{\rm TG}}{3} a^3 e^{-a^2/2} \right)\mathbf{I}_3 \label{eq:L_TG_3D}
\end{align}

Finally, the CRLLB is given by
\begin{align}
\mathbf{P}({\bf x})& \geq  \mathbf{L}({\bf x})^{\top}\mathbf{J}({\bf x})^{-1}\mathbf{L}({\bf x})\nonumber\\
& = \left(1 -  \frac{4\pi k^3_{\rm TG}}{3} a^3 e^{-a^2/2} \right)\left( \mathbf{H}^{\top} \mathbf{H}\right)^{-1}\label{eq:CRLLB_TG_b}
\end{align}
which is the same as its covariance \eqref{eq:cov_TG_3D}. This is, again, because its LF belongs to the exponential family and it satisfies \textit{generalized collinearity condition} \eqref{collinear}.

\section{Estimation of the support parameters of a uniform distribution}\label{sec:uniform_example}
The aforementioned examples are all about additive noise models~\eqref{eq:model}. Next, we investigate another special type of LF~\cite{ye2020estimation}. Consider $N$ i.i.d. measurements ${\bf Z} = [z_1,\ldots, z_N]^\top$ drawn from the uniform distribution ${\cal U}(x_1,x_2)$, with unknown parameter vector ${\bf x} \triangleq [x_1, x_2]^\top$ to be estimated. The LF is then given by
\begin{align}
   p({\bf Z}|{\bf x}) = \prod_{n=1}^N p(z_n|x_1,x_2) = \frac{1}{(x_2-x_1)^n}, z_n\in [x_1,x_2] \label{eq:unif_LF}
\end{align}
The MLEs for $x_1$ and $x_2$ are given by
\begin{align}
   z_m = {\rm min}({\bf Z})  \quad z_M = {\rm max}({\bf Z})
\end{align}
which are not unbiased estimators of ${\bf x}$. To derive the unbiased estimators, we will first obtain the pdfs of $z_m$ and $z_M$, which are
\begin{align}
    p(z_m|{\bf x})& = \frac{n(x_2 - z_m)^{n-1}}{(x_2-x_1)^n}\quad x_1\leq z_m \leq x_2\\
      p(z_M|{\bf x}) &= \frac{n(z_M - x_1)^{n-1}}{(x_2-x_1)^n}\quad x_1\leq z_M \leq x_2
\end{align}
whose expected values are given by
\begin{align}
    \mathbb{E}[z_m|{\bf x}] & = \frac{nx_1 + x_2}{n+1}\\
    \mathbb{E}[z_M|{\bf x}] & = \frac{nx_2 + x_1}{n+1}
\end{align}
Hence, one can obtain the unbiased estimators for ${\bf x}$ as
\begin{align}
    \hat{x}_1 ({\bf z}) & = \frac{nz_m - z_M}{n-1}\\
    \hat{x}_2 ({\bf z}) & = \frac{nz_M - z_m}{n-1}
\end{align}
where ${\bf z}  \triangleq [z_m, z_M]^\top$ is the sufficient statistic containing all relevant information about ${\bf x}$. Hereafter, we will use ${\bf z}$ as the transformed measurements, yielding the LF
\begin{align}
    p({\bf z}| {\bf x}) = n(n-1)\frac{(z_M - z_m)^{n-2}}{(x_2 - x_1)^n} \label{eq: LF_uniform}
\end{align}
which has parameter dependent support 
\begin{align}
 {\cal S}({\bf x}) = \{{\bf x}|x_1\leq z_m \leq z_M \leq x_2 \} \label{eq:support_uniform}
\end{align}
Having available~\eqref{eq: LF_uniform} and~\eqref{eq:support_uniform}, the covariance matrix for $\hat{x}$ is given by
\begin{align}
   {\bf P}({\bf x}) &= \int_{{\cal S}({\bf x})} p({\bf z}|{\bf x})(\hat{\bf x}({\bf z}) - {\bf x})(\hat{\bf x}({\bf z}) - {\bf x})^\top d{\bf z} \nonumber\\
   & = \begin{bmatrix}
       P_{11}({\bf x}) & P_{12}({\bf x})\\
       P_{21}({\bf x}) & P_{22} ({\bf x})
   \end{bmatrix}
\end{align}
where 
\begin{align}
    P_{11}({\bf x}) & = P_{22}({\bf x}) = \frac{n(x_2 - x_1)^2}{(n-1)(n+1)(n+2)}\\
    P_{12}({\bf x}) & = P_{21}({\bf x}) = \frac{-(x_2 - x_1)^2}{(n-1)(n+1)(n+2)}
\end{align}

\begin{figure}[!t]
    \centering
    \includegraphics[width=0.35\textwidth]{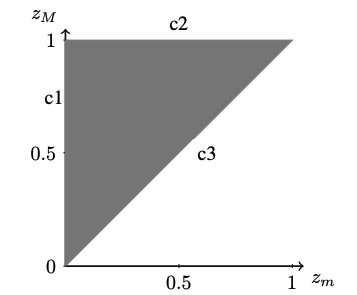}
    \caption{Support ${\cal S}({\bf x})$ of the uniform LF for $x_1=0$ and $x_2=1$.}
    \label{fig:uniform_support}
\end{figure}

Based on the LF~\eqref{eq:unif_LF}, the FIM is calculated as
\begin{align}
    {\bf J}({\bf x}) &= \mathbb{E}[\nabla_{\bf x} \ln p({\bf Z}|{\bf x}) \nabla_{\bf x} \ln p({\bf Z}|{\bf x})^\top] \nonumber\\
    & = \begin{bmatrix}
        \frac{n^2}{(x_2 - x_1)^2} & -\frac{n^2}{(x_2 - x_1)^2}\\
        -\frac{n^2}{(x_2 - x_1)^2} & \frac{n^2}{(x_2 - x_1)^2}
    \end{bmatrix}
\end{align}
which, apparently, is singular and not invertible. Thus, there exists no meaningful CRLB or CRLLB in this case. To get around this hurdle and obtain a valid variance bound, we will approximate the uniform LF with a truncated Gaussian (TG) LF with variance large enough so as to yield a flat distribution over its support. To proceed, the TG-approximated LF is given by (with approximation indicated by an apostrophe)
\begin{align}
    p'(z_n|{\bf x}) = \frac{{\cal N}(z_n; \mu, \sigma^2)}{\Phi(x_2; \mu, \sigma^2) - \Phi(x_1; \mu, \sigma^2)}, x_1 \leq z_n\leq x_2
\end{align}
where ${\cal N}(z_n; \mu, \sigma^2)$ is a Gaussian pdf with mean $\mu = (x_1 +x_2)/2$ and variance $\sigma^2$, and $\Phi(x; \mu, \sigma^2) = \int_{-\infty}^x {\cal N}(z; \mu, \sigma^2) dz$ is a Gaussian cumulative density function (cdf).

The LF for the sufficient statistic vector ${\bf z}$ is then given by
\begin{align}
   p'({\bf z}| {\bf x}) & =  n(n-1){\cal N}(z_m; \mu, \sigma^2){\cal N}(z_M; \mu, \sigma^2) \nonumber\\
  &\  \quad\frac{[\Phi(z_M; \mu, \sigma^2) - \Phi(z_m; \mu, \sigma^2)]^{n-2}}{[\Phi(x_2; \mu, \sigma^2) - \Phi(x_1; \mu, \sigma^2)]^n} \label{eq:unif_LF_TG}
\end{align}
Evaluating the FIM with the new LF yields 
\begin{align}
   {\bf J}'({\bf x})= \begin{bmatrix}
       J'_{11}({\bf x}) & J'_{12}({\bf x}) \\
       J'_{21}({\bf x}) & J'_{22}({\bf x})
    \end{bmatrix} \label{eq:FIM_TG}
\end{align}
where
\begin{align}
 & J'_{11}({\bf x}) = J'_{22}({\bf x}) = \frac{n^2{\cal N}(\frac{x_1-x_2}{2};0, \sigma^2)^2}{[1-2\Phi(\frac{x_1-x_2}{2};0, \sigma^2)]^2} + \frac{n}{4\sigma^2} \\
 & J'_{12}({\bf x}) = J'_{21}({\bf x}) =  -\frac{n^2{\cal N}(\frac{x_1-x_2}{2};0, \sigma^2)^2}{[1-2\Phi(\frac{x_1-x_2}{2};0, \sigma^2)]^2} + \frac{n}{4\sigma^2}
\end{align}
It shall be easy to verify that 
\begin{align}
 \lim_{\sigma^2\rightarrow \infty}  {\bf J}'({\bf x}) = {\bf J}({\bf x}).
\end{align}
When $\sigma^2$ is finite, the FIM~\eqref{eq:FIM_TG} is invertible, which results in the existence of the CRLB. Nevertheless, the fact that the support is parameter-dependent renders the CRLB not a valid bound. Next, we will evaluate the Leibniz term to obtain the CRLLB using the formula~\eqref{eq:D_L_2}. As shown in Fig.~\ref{fig:uniform_support}, there are three segments ($c1$, $c2$, and $c3$) along the boundary. The corresponding $\nabla_{\bf x} {\bf z}^\top$ and ${\bf n}$ are respectively given by
\begin{align}
 \nabla_{\bf x} {\bf z}_1^\top & = \begin{bmatrix}
     1 & \frac{x_2 - z_M}{x_2 - x_1}\\ 0 & \frac{z_M - x_1}{x_2 - x_1}
 \end{bmatrix}  \quad 
 & {\bf n}^{(1)}  = \begin{bmatrix}
     -1 \\ 0
 \end{bmatrix} \nonumber\\
  \nabla_{\bf x} {\bf z}_2^\top & = \begin{bmatrix}
     \frac{x_2 - z_m}{x_2 - x_1} &  0 \\ \frac{z_m - x_1}{x_2 - x_1} & 1 
 \end{bmatrix}  \quad 
& {\bf n}^{(2)}  = \begin{bmatrix}
     0 \\ 1
 \end{bmatrix} \nonumber\\
\nabla_{\bf x} {\bf z}_3^\top & = \begin{bmatrix}
     \frac{x_2 - z_m}{x_2 - x_1} &  \frac{x_2 - z_M}{x_2 - x_1} \\ \frac{z_m - x_1}{x_2 - x_1} &  \frac{z_M - x_1}{x_2 - x_1}
 \end{bmatrix}  \quad 
& {\bf n}^{(3)}  = \begin{bmatrix}
     \frac{\sqrt{2}}{2} \\ -\frac{\sqrt{2}}{2}
 \end{bmatrix} 
\end{align}
Using $p'({\bf z}|{\bf x})$~\eqref{eq:unif_LF_TG}, the Leibniz term is then given by
\begin{align}
   & {\bf D}'_{L_1} ({\bf x}) = \int_{c1}  \nabla_{\bf x} {\bf z}_1^\top {\bf n}^{(1)} p'({\bf z}|{\bf x})(\hat{\bf x}({\bf z})-{\bf x})^\top dc \nonumber\\
   & = \int_{x_1}^{x_2} \begin{bmatrix}
       \frac{z_M - x_1}{n-1} & 0\\ 0 & 0
   \end{bmatrix} n(n-1){\cal N}(x_1; \mu, \sigma^2){\cal N}(z_M; \mu, \sigma^2) \nonumber\\
  &\  \quad\frac{[\Phi(z_M; \mu, \sigma^2) - \Phi(x_1; \mu, \sigma^2)]^{n-2}}{[\Phi(x_2; \mu, \sigma^2) - \Phi(x_1; \mu, \sigma^2)]^n} d z_M
\end{align}
Integrating out terms related to $z_M$, we have (the arguments $\mu$ and $\sigma^2$ in ${\cal N}$ and $\Phi$ are dropped for notational brevity)
\begin{align}
   & \int_{x_1}^{x_2} (z_M - x_1){\cal N}(z_M)[\Phi(z_M) - \Phi(x_1)]^{n-2} d z_M \nonumber\\
   & = \int_{x_1}^{x_2} (z_M - x_1) d \frac{[\Phi(z_M) - \Phi(x_1)]^{n-1}}{n-1} \nonumber\\
   & = (x_2 - x_1)\frac{[\Phi(x_2) - \Phi(x_1)]^{n-1}}{n-1} \nonumber\\
   & \quad \ - \int_{x_1}^{x_2} \frac{[\Phi(z_M) - \Phi(x_1)]^{n-1}}{n-1} d z_M
\end{align}
Thus,
\begin{align}
  & {\bf D}'_{L_1} ({\bf x}) = \begin{bmatrix}
       1 & 0\\ 0 & 0
   \end{bmatrix} \Bigg(\frac{n{\cal N}(x_1; \mu, \sigma^2)(x_2 - x_1)}{(n-1)[\Phi(x_2; \mu, \sigma^2) - \Phi(x_1; \mu, \sigma^2)]}\nonumber\\
   & - \frac{n{\cal N}(x_1; \mu, \sigma^2)\int_{x_1}^{x_2} [\Phi(z_M; \mu, \sigma^2) - \Phi(x_1; \mu, \sigma^2)]^{n-1} d z_M}{(n-1)[\Phi(x_2; \mu, \sigma^2) - \Phi(x_1; \mu, \sigma^2)]^n}\Bigg)
\end{align}
Similarly, we have
\begin{align}
  & {\bf D}'_{L_2} ({\bf x}) = \begin{bmatrix}
       0 & 0\\ 0 & 1
   \end{bmatrix} \Bigg(\frac{n{\cal N}(x_2; \mu, \sigma^2)(x_2 - x_1)}{(n-1)[\Phi(x_2; \mu, \sigma^2) - \Phi(x_1; \mu, \sigma^2)]}\nonumber\\
   & - \frac{n{\cal N}(x_2; \mu, \sigma^2)\int_{x_1}^{x_2} [\Phi(x_2; \mu, \sigma^2) - \Phi(z_m; \mu, \sigma^2)]^{n-1} d z_m}{(n-1)[\Phi(x_2; \mu, \sigma^2) - \Phi(x_1; \mu, \sigma^2)]^n} \Bigg)
\end{align}
Thus, ${\bf L}'({\bf x})$ becomes
\begin{align}
  {\bf L}'({\bf x}) = {\bf I}_2 -{\bf D}'_{L_1} ({\bf x})-{\bf D}'_{L_2} ({\bf x}) 
\end{align}
Apparently, as $\sigma$ approaches infinity, the following two hold
\begin{align}
  & \lim_{\sigma\rightarrow\infty} \Phi(x_2; \mu, \sigma^2) - \Phi(x_1; \mu, \sigma^2) = \frac{x_2 - x_1}{\sqrt{2\pi}\sigma}\\
  & \lim_{\sigma\rightarrow\infty}{\cal N}(x;\mu, \sigma^2)  = \frac{1}{\sqrt{2\pi}\sigma}
\end{align}
which will result in
\begin{align}
&{\bf D}_{L_1} ({\bf x}) = \lim_{\sigma\rightarrow\infty} {\bf D}'_{L_1} ({\bf x})
=\begin{bmatrix}
       1 & 0\\ 0 & 0
   \end{bmatrix}\\
   &{\bf D}_{L_2} ({\bf x})= \lim_{\sigma\rightarrow\infty} {\bf D}'_{L_2} ({\bf x})
 = \begin{bmatrix}
       0 & 0\\ 0 & 1
   \end{bmatrix}
\end{align}
For the uniform LF, the following holds
\begin{align}
  {\bf L}({\bf x}) = {\bf I}_2 -{\bf D}_{L_1} ({\bf x})-{\bf D}_{L_2} ({\bf x}) = {\bf 0}  
\end{align}
which yields the corresponding CRLLB being ${\bf 0}$, as derived in~\cite{ye2020estimation}. In spite of being ${\bf 0}$, the CRLLB is still a valid bound.

\section{Conclusion and Discussion}\label{sec:conclusion}
This paper relies on the general framework of the multidimensional  CRLLB to unify existing results on CRLB and CRLLB, including the scenarios with scalar parameters and LFs with parameter-dependent support. 
Notably, the regularity condition for the CRLLB to hold with equality is the \textit{generalized collinearity condition} \eqref{collinear} between the gradient of LLF w.r.t. ${\bf x}$ and the estimation error, in which case the LF has to belong to the exponential family. To demonstrate the usefulness of the CRLLB, some typical examples of LFs have been surveyed. For the case of the RFC LF inside a circle, the CRLLB provides a valid bound but it is not achievable as the LF does not satisfy the \textit{generalized collinearity condition}  (it does not belong to the exponential family). If the measurement noise has a TG distribution inside an $n$-sphere, the variance of the unbiased estimator equals the CRLLB because the \textit{generalized collinearity condition} is satisfied. For measurement noise with a truncated Laplace distribution within a circle, although the LF belongs to the exponential family, the variance is larger than the CRLLB for that the \textit{generalized collinearity condition}  is not satisfied. This demonstrates that the LF being exponential is only an necessary condition for the CRLLB to hold with equality but not a sufficient condition. Interestingly, the uniform distribution inside a sphere, which can either be obtained by the RFC pdf with $\beta = 0$ or be approximated by the TG pdf as $\sigma\rightarrow\infty$. Depending on the form of the approximated LF, the CRLLB and the claim of statistical efficiency varies.
In addition, we provide an example of an observation model with non-identity measurement matrix with $3$-D TG distribution for the measurement noise. When the dimension of the measurement is larger than the dimension of the parameter, the MLE is an efficient estimator as the \textit{generalized collinearity condition} holds. Lastly, we revisit the example of the estimation of the support parameters of a uniform distribution with a valid CRLLB derived.

A promising application domain of the CRLLB is in quantum metrology or quantum parameter estimation~\cite{giovannetti2011advances}, where the parameter of interest takes values that change the rank of the quantum statistical model (similar to the case with parameter-dependent support here)~\cite{seveso2019discontinuity, ciaglia2023categorical, ciaglia2020differential}. How to adapt the CRLLB to the quantum counterparts will be an interesting direction to pursue.
\appendices
\section{Generalization of Leibniz rule for the gradient of a multidimensional integration \label{sec:Leibniz}\protect\footnote{based on \cite{Flanders} pp 618--624}}
\subsection{Differentiation of a 2-dimensional integral (the plane formula)}
\begin{figure*}[!t]
\centering
\subfloat[]{\includegraphics[width=0.4\textwidth]{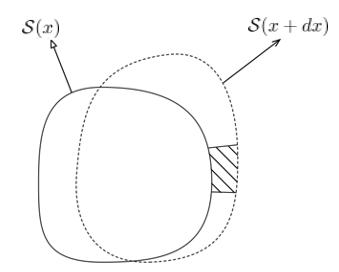}%
\label{fig:S_t}}
\hfil
\subfloat[]{\includegraphics[width=0.25\textwidth]{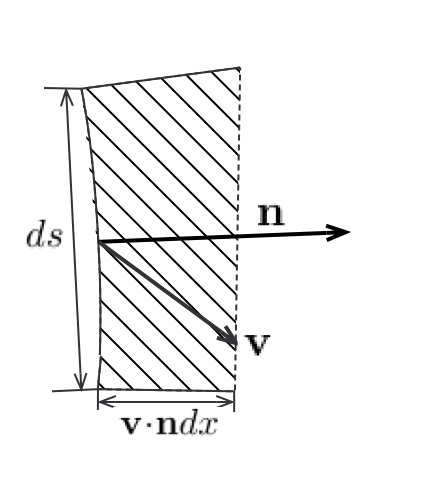}%
\label{fig:deltaDt}}
\caption{Illustration of two successive domains}
\label{fig:deltaS_t}
\end{figure*}
Suppose we have a domain $\mathcal{S}(x)$ in the $z_1$-$z_2$ plane that depends on an additional variable $x$. We are given a function $F(z_1,z_2,x)$ and the problem is to find
\begin{equation}\label{eq:2D_integration}
\frac{d}{dx}\int\!\!\!\int_{\mathcal{S}(x)}F(z_1,z_2,x)dz_1 dz_2\Bigr\rvert_{x= x_0}
\end{equation}
The first move should be separation of the boundary variation from integrand variation, which is given by

\begin{align}
&\frac{d}{dx}\int\!\!\!\int_{\mathcal{S}(x)}F(z_1,z_2,x)dz_1 dz_2\Bigr\rvert_{x= x_0}\nonumber \\
 = & \int\!\!\!\int_{\mathcal{S}(x)}\frac{\partial F}{\partial x}\Bigr\rvert_{x= x_0} dz_1 dz_2 + \nonumber \\
 & \frac{1}{dx}\left(\int\!\!\!\int_{\mathcal{S}(x+dx)} - \int\!\!\!\int_{\mathcal{S}(x)}\right) F(z_1,z_2,x) dz_1 dz_2\Bigr\rvert_{x= x_0}
\end{align}

The essence of the problem is to find the second term, i.e. the Leibniz term, in the above equation. 
Two successive domains $\mathcal{S}(x)$ and $\mathcal{S}(x+dx)$ are illustrated in Fig.\ \ref{fig:S_t}. Let 
\begin{equation}
\mathbf{v} = \left[\frac{dz_1}{dx}\quad \frac{dz_2}{dx}\right]^{\top} = [v_1\quad v_2]^{\top}
\end{equation}
denote the velocity vector (if we notionally regard $x$ as time) at a boundary point $(z_1, z_2)$ of $\mathcal{S}(x)$ and $\mathbf{n}$ denote the outward unit normal. The infinitesimal displacement of the boundary (contour) $ds$ is the shaded region in Fig.\ \ref{fig:deltaS_t} and its area is given by the projection of the velocity vector $\mathbf{v}$ multiplied by $dx$ on the normal to the boundary $\mathbf{n}$. The Leibniz term is then
\begin{align}\label{eq:Leibniz_2D_initial}
&\frac{1}{dx}\left(\int\!\!\!\int_{\mathcal{S}(x+dx)} - \int\!\!\!\int_{\mathcal{S}(x)}\right) F(z_1,z_2,x) dz_1dz_2\nonumber\\
 &\approx \int_{\delta \mathcal{S}(x)} F(z_1, z_2) \mathbf{v}\cdot \mathbf{n}ds
\end{align}
where $\delta$ denotes the boundary, $\mathbf{v}\cdot \mathbf{n}$ is the dot (inner) product of the vectors $\mathbf{v}$ and $\mathbf{n}$, and $ds$ is the infinitesimal contour of the boundary. 

Rotating the unit tangent $[dz_1/ds\quad dz_2/ds]^{\top}$ to $\delta \mathcal{S}_x$ counter clockwise through $90^{\circ}$ yields 
\begin{equation}
\mathbf{n} = [dz_1/ds\quad -dz_2/ds]^{\top}
\end{equation}
Hence, the Leibniz term in the plane is 
\begin{align}
&\frac{1}{dx}\left(\int\!\!\!\int_{\mathcal{S}(x+dx)} - \int\!\!\!\int_{\mathcal{S}(x)}\right) F(z_1, z_2, x) dz_1dz_2 \nonumber\\
& = \int_{\delta \mathcal{S}(x)} F(z_1, z_2)(v_1 dz_2 - v_2 dz_1)\label{eq:Leibniz_2D}
\end{align}

\subsection{Differentiation of a 3-dimensional integral}
Suppose we are given a function $F(z_1, z_2, z_3, x)$, and the problem is to find the differentiation w.r.t. $x$ under a 3-dimensional integral of $F$ w.r.t. $z_1$, $z_2$ and $z_3$ with the integration space dependent on $x$, i.e.,
\begin{equation}\label{eq:3D_integration}
\frac{d}{dx}\int\!\!\!\int_{\mathcal{S}(x)}F(z_1,z_2, z_3, x)dz_1 dz_2 dz_3\Bigr\rvert_{x= x_0}
\end{equation}
As with the previous subsection, the first step is the separation of the variations of the integrand with the boundary of the integration space, which is a surface. 
At any point along the boundary $\delta \mathcal{S}(x)$, let  $ds$ represent infinitesimal area around that point, $\mathbf{n}$ denote the unit vector normal to the plane and  $\mathbf{v}$ represent the 3-dimensional velocity vector (if regard $x$ as time) as
\begin{equation}
\mathbf{v} = \left[\frac{dz_1}{dx}\quad \frac{dz_2}{dx}\quad\frac{dz_3}{dx}\right]^{\top} = \left[v_1\quad v_2 \quad v_3\right]^{\top}
\end{equation}
The variation of the integration volume is multiplying $ds$ by the projection of $\mathbf{v}$ on $\mathbf{n}$. Therefore,
\begin{subequations}
\begin{align}
&\frac{1}{dx}\left(\int\!\!\!\int_{\mathcal{S}(x+dx)} - \int\!\!\!\int_{\mathcal{S}(x)}\right) F(z_1, z_2, z_3, x) dz_1dz_2 dz_3 \nonumber\\
&= \int_{\delta \mathcal{S}(x)} F \mathbf{v}\cdot \mathbf{n}ds \label{eq:Leibniz_3D_b}\\
& = \int_{\delta \mathcal{S}(x)} \left( Fv_1 dz_2 dz_3 + Fv_2 dz_3 dz_1  + Fv_3 dz_1 dz_2 \right)\label{eq:Leibniz_3D_c}
\end{align}
\end{subequations}
The derivation from \eqref{eq:Leibniz_3D_b} to \eqref{eq:Leibniz_3D_c} follows from the following identity:
\begin{equation}
\mathbf{n}ds = \left[dz_2 dz_3 \quad dz_3 dz_1  \quad dz_1 dz_2 \right]^{\top}
\end{equation}
\subsection{Generalization of differentiation of a multidimensional integral}\label{sec:general_Leibniz}
Differentiation under 2-dimensional \eqref{eq:2D_integration} and 3-dimensional \eqref{eq:3D_integration} integrals can be obtained with the Euclidean structure of space. But when it comes to higher-dimensional integral, it is not feasible to do so. To generalize the Leibniz rule for differentiation of a multidimensional integral, we should get acquainted with the following mathematical terminologies
\begin{enumerate}
\item \textbf{Vector field}: a vector field is an assignment of a vector to each point in a subset of a space.
\item \textbf{Exterior product}: the exterior product of a set of vectors gives the volume of the subspace the vectors define and also the orientation of the subspace.
\item \textbf{Differential form}: differential forms provide a unified approach to defining integrands over curves, surfaces, volumes, and higher-dimensional manifolds. For instance, $f(x) dx$ is an example of $1$-form that can be integrated over interval $[a, b]$; $f(x, y, z) dx \wedge dy + g(x, y, z) dx \wedge dz + h(x, y, z) dy \wedge dz$ is a $2$-form that can be integrated over a surface.
\item \textbf{Interior product}: the interior product compresses a $p$-form to a $(p-1)$-form in a smooth manifold.
\end{enumerate}

Suppose we want to obtain the following differentiation with respect to $x$ of a $n_{\mathbf{z}}$-dimensional integral,
\begin{equation}
\frac{d}{dx}\int_{\mathcal{S}(x)}F(\mathbf{z},x) d\mathbf{z}
\end{equation}
where $\mathcal{S}(x)$ is the integration space, 
\begin{equation}
\mathbf{z} = [z_1\quad ... \quad z_{n_{\mathbf{z}}}]^{\top}
\end{equation}
and $d\mathbf{z}$ is a $n_z$-form, given by 
\begin{equation}
d\mathbf{z} = dz_{1}\wedge dz_{2}\wedge...\wedge dz_{n_{\mathbf{z}}}
\end{equation}
in which $\wedge$ represents the exterior product. 

Let $\omega = F(\mathbf{z},x) d\mathbf{z}$, the general Leibniz term for differentiation under $n_{\mathbf{z}}$-dimensional integral is given by \cite{Flanders}
\begin{equation}\label{eq:general Leibniz}
 \delta_x \left( \frac{d}{dx}\int_{\mathcal{S}(x)}F(\mathbf{z},x) d\mathbf{z}\right) = \int_{\delta{\mathcal{S}(x)}} \left(\mathbf{v}\cdot \nabla_{\mathbf{z}}\right)\lrcorner\, \omega
\end{equation}
where $\delta_x$ means contribution from the variation of the boundary w.r.t. $x$, $\mathbf{v}$ as defined in the case of $n$-dimensional integral is the differentiation of the point on the boundary of the integration space w.r.t. $x$, i.e.,
\begin{align}
\mathbf{v} \triangleq \frac{d \mathbf{z}}{d x} = \left[\frac{dz_1}{dx}\quad \frac{dz_2}{dx}\quad .... \quad \frac{dz_{n_{\mathbf{z}}}}{dx}\right]^{\top}
\end{align}
and $\nabla_{\mathbf{z}}$ is the operator
\begin{align}
\nabla_{\mathbf{z}} = \left[ \frac{\partial}{\partial z_1} \quad \frac{\partial}{\partial z_2}\quad ... \quad \frac{\partial}{\partial z_{n_{\mathbf{z}}}} \right]
^{\top}
\end{align}
$\mathbf{v}\cdot \nabla_{\mathbf{z}}$ defines a vector field.
 The details of the proof of the formula are not provided here as it is lengthy and rather technical. Interested readers are encouraged to read \cite{Flanders}.

In \eqref{eq:general Leibniz}, $\lrcorner$ is the interior product operator, which contracts a $n_{\mathbf{z}}$-form to a $(n_{\mathbf{z}} - 1)$-form within a vector field. The interior product of $\mathbf{v}$ and a decomposable $n_z$-form $\alpha^1\wedge\alpha^2\wedge...\wedge\alpha^{n_{\mathbf{z}}}$ is defined by
\begin{align}\label{eq:interior product}
&\left(\mathbf{v}\cdot \nabla_{\mathbf{z}}\right)\lrcorner(\alpha^1\wedge\alpha^2\wedge...\wedge\alpha^{n_{\mathbf{z}}}) =\nonumber\\
& \sum_i (-1)^{i-1}\left<\mathbf{v}\cdot \nabla_{\mathbf{z}},\alpha^i\right>\alpha^1\wedge...\wedge\alpha^{i-1}\wedge\alpha^{i+1}\wedge...\wedge\alpha^{n_{\mathbf{z}}}
\end{align}
where $\alpha^i$ is one-form, denoted as
\begin{align}
\alpha^i = \sum_{j = 1}^{n_{\mathbf{z}}} a_j^i dz_j
\end{align}
and 
$\left<\mathbf{v}\cdot \nabla_{\mathbf{z}},\alpha^i\right>$ determines the effect of the one-form $\alpha^i$ on a vector field $\mathbf{v}$, which is
\begin{equation}
<\mathbf{v}\cdot \nabla_{\mathbf{z}},\alpha^i> = \left<\sum_{j = 1}^{n_{\mathbf{z}}} v_{j}\frac{\partial}{\partial z_j}, \sum_{j = 1}^{n_{\mathbf{z}}}a_j^i dz_j\right> = \sum_{j = 1}^{n_{\mathbf{z}}} v_j a_j^i
\end{equation}
The Leibniz terms given by \eqref{eq:Leibniz_2D} and \eqref{eq:Leibniz_3D_c} are two particular cases of \eqref{eq:general Leibniz}. 

For differentiation of a 2-dimensional integral in \eqref{eq:2D_integration},
\begin{equation}
 \omega = F(z_1, z_2, x)dz_1\wedge dz_2
\end{equation}
and the Leibniz term according to \eqref{eq:general Leibniz} is
\begin{align}
&\int_{\delta \mathcal{S}(x)} \left(\mathbf{v}\cdot \nabla_{\mathbf{z}}\right)\lrcorner\omega  = \int_{\delta\mathcal{S}(x)} \left(\mathbf{v}\cdot \nabla_{\mathbf{z}}\right)\lrcorner\left[F(z_1, z_2, x)dz_1\wedge dz_2 \right]\nonumber\\
& = \int_{\delta \mathcal{S}(x)} F\Bigg[ \left<v_1\frac{\partial}{\partial z_1} + v_2\frac{\partial}{\partial z_2}, dz_1\right>dz_2 \nonumber\\
& \quad- \left<v_1\frac{\partial}{\partial z_1} + v_2\frac{\partial}{\partial z_2},dz_2\right>dz_1\Bigg]\nonumber \\
& = \int_{\delta  \mathcal{S}(x)} F\left(v_1 dz_2-v_2 dz_1\right)\label{eq:2D_general_c}
\end{align}
which is the same as \eqref{eq:Leibniz_2D}. 

For differentiation of a 3-dimensional integral in \eqref{eq:3D_integration},
\begin{equation}
 \omega = F(z_1, z_2, z_3, x)dz_1\wedge dz_2 \wedge dz_3
\end{equation}
and the Leibniz term according to \eqref{eq:general Leibniz} is
\begin{align}
&\int_{\delta \mathcal{S}(x)} \left(\mathbf{v}\cdot \nabla_{\mathbf{z}}\right)\lrcorner\omega \nonumber\\
 =& \int_{\delta\mathcal{S}(x)} \left(\mathbf{v}\cdot \nabla_{\mathbf{z}}\right)\lrcorner\left[F(z_1, z_2, z_3, x)dz_1\wedge dz_2  \wedge dz_3\right]\nonumber\\
= &\int_{\delta \mathcal{S}(x)} F\Bigg[ \left<v_1\frac{\partial}{\partial z_1} + v_2\frac{\partial}{\partial z_2}+ v_3\frac{\partial}{\partial z_3}, dz_1\right>dz_2\wedge dz_3 \nonumber\\
&- \left<v_1\frac{\partial}{\partial z_1} + v_2\frac{\partial}{\partial z_2}+ v_3\frac{\partial}{\partial z_3},dz_2\right>dz_1\wedge dz_3 \nonumber\\
&+ \left<v_1\frac{\partial}{\partial z_1} + v_2\frac{\partial}{\partial z_2}+ v_3\frac{\partial}{\partial z_3},dz_3\right>dz_1\wedge dz_2\Bigg]\nonumber\\
 =& \int_{\delta \mathcal{S}(x)} \left( Fv_1 dz_2 \wedge dz_3 - Fv_2 dz_1\wedge dz_3  + Fv_3 dz_1\wedge dz_2 \right)\nonumber\\
 = & \int_{\delta \mathcal{S}(x)} \left( Fv_1 dz_2 dz_3 + Fv_2 dz_3 dz_1  + Fv_3 dz_1 dz_2 \right)\label{eq:3D_general_d}
\end{align}
which is the same as \eqref{eq:Leibniz_3D_c}. 

\subsection{Gradient of a multidimensional integral}
Suppose we have a function $F(\mathbf{z}, \mathbf{x})$ dependent on $n_{\mathbf{z}}$-dimensional vector $\mathbf{z} = [z_1\quad \cdots\quad z_{n_{\mathbf{z}}}]^{T}$ and $n_{\mathbf{x}}$-dimensional vector $\mathbf{x} = [x_1\quad \cdots\quad x_{n_{\mathbf{x}}}]^{T}$ and we are interested in
\begin{align}\label{eq:general_gradient_integral}
\nabla_{\mathbf{x}} \int_{\mathcal{S}(\mathbf{x})} F(\mathbf{z}, \mathbf{x}) d\mathbf{z}
= \begin{bmatrix}
\frac{\partial}{\partial x_1}\\
\vdots \\
\frac{\partial}{\partial x_{n_{\mathbf{x}}}}
\end{bmatrix}\int_{\mathcal{S}(\mathbf{x})} F(\mathbf{z}, \mathbf{x}) d\mathbf{z}
\end{align}
Let $\omega = F(\mathbf{z}, \mathbf{x})d\mathbf{z}$, the Leibniz term (the contribution from the integration space) of the above gradient w.r.t. $\mathbf{x}$ under integration w.r.t. $\mathbf{z}$ is given according to \eqref{eq:general Leibniz} by
\begin{align}\label{eq:Leibniz_gradient}
\delta_{\mathbf{x}} \left(\nabla_{\mathbf{x}} \int_{\mathcal{S}(\mathbf{x})} F(\mathbf{z}, \mathbf{x}) d\mathbf{z}\right) = \int_{\delta\mathcal{S}(\mathbf{x})}\begin{bmatrix}
\left(\mathbf{v}_1\cdot \nabla_{\mathbf{z}}\right) \lrcorner \omega\\
\vdots\\
\left(\mathbf{v}_q\cdot \nabla_{\mathbf{z}}\right) \lrcorner \omega
\end{bmatrix}
\end{align}
where $\delta_{\mathbf{x}}$ means contribution from the variation of the boundary w.r.t. $\mathbf{x}$ and $\mathbf{v}_j$, $j = 1, ..., n_{\mathbf{x}}$, is
\begin{align}
\mathbf{v}_j = \frac{\partial}{\partial x_j} \mathbf{z} = \left[\frac{\partial z_1}{\partial x_j}\quad\frac{\partial z_2}{\partial x_j}\quad .... \quad \frac{\partial z_{n_{\mathbf{z}}}}{\partial x_j}\right]^{\top}
\end{align}
And \eqref{eq:Leibniz_gradient} can be rewritten using matrix operation as
\begin{subequations}
\begin{align}
&\delta_{\mathbf{x}} \left( \nabla_{\mathbf{x}} \int_{\mathcal{S}(\mathbf{x})} F(\mathbf{z}, \mathbf{x}) d\mathbf{z}\right) = \int_{\delta\mathcal{S}(\mathbf{x})} \left(\nabla_{\mathbf{x}} \mathbf{z}^{T} \nabla_{\mathbf{z}}\right) \lrcorner \omega \\
& = \int_{\delta\mathcal{S}(\mathbf{x})}\left(\begin{bmatrix}
\frac{\partial z_1}{\partial x_1} & \cdots & \frac{\partial z_{n_{\mathbf{z}}}}{\partial x_1}\\
\vdots & \ddots & \vdots\\
\frac{\partial z_1}{\partial x_{n_{\mathbf{x}}}} & \cdots & \frac{\partial z_{n_{\mathbf{z}}}}{\partial x_{n_{\mathbf{x}}}}
\end{bmatrix}
\begin{bmatrix}
\frac{\partial}{\partial z_1} \\ \vdots \\ \frac{\partial}{\partial z_{n_{\mathbf{z}}}} 
\end{bmatrix} \right)\lrcorner \omega \label{eq:Leibniz_gradient_matrix_b}
\end{align}
\end{subequations}
where \eqref{eq:Leibniz_gradient_matrix_b} is the same as \eqref{eq:Leibniz_gradient}.

When we have a $1\times n_{\mathbf{F}}$ function vector
\begin{equation}
\mathbf{F}(\mathbf{z}, \mathbf{x}) = \left[F_1 (\mathbf{z}, \mathbf{x})\quad \cdots \quad F_{n_{\mathbf{F}}} (\mathbf{z}, \mathbf{x})\right]
\end{equation}
and let
\begin{equation}
\boldsymbol{\omega} = \left[F_1 (\mathbf{z}, \mathbf{x})d\mathbf{z} \quad \cdots \quad F_{n_{\mathbf{F}}} (\mathbf{z}, \mathbf{x})d\mathbf{z}  \right]
\end{equation}
Then the gradient w.r.t. $\mathbf{x}$ under the multidimensional integral for a vector function is given by
\begin{align}
&\delta_{\mathbf{x}} \left(\nabla_{\mathbf{x}} \int_{\mathcal{S}(\mathbf{x})} \mathbf{F}(\mathbf{z}, \mathbf{x}) d\mathbf{z}\right)\nonumber\\
& = \int_{\delta\mathcal{S}(\mathbf{x})} 
\left[ \left(\nabla_{\mathbf{x}} \mathbf{z}^{T} \nabla_{\mathbf{z}}\right) \lrcorner \omega_1\quad \cdots \quad \left(\nabla_{\mathbf{x}} \mathbf{z}^{T} \nabla_{\mathbf{z}}\right) \lrcorner \omega_{n_{\mathbf{F}}} \right]\nonumber\\
&\triangleq  \int_{\delta\mathcal{S}(\mathbf{x})} \left(\nabla_{\mathbf{x}} \mathbf{z}^{T} \nabla_{\mathbf{z}}\right) \lrcorner \boldsymbol{\omega} \label{eq:Leibniz_gradient_matrix_a}
\end{align}

When the dimension of $\mathbf{z}$ is $2$ or $3$, we can have a simplified form (without the expression of the interior product) as 
\begin{equation}\label{eq:Leibniz_2D_3D}
\delta_{\mathbf{x}} \left( \nabla_{\mathbf{x}} \int_{\mathcal{S}(\mathbf{x})} \mathbf{F}(\mathbf{z}, \mathbf{x}) d\mathbf{z}\right)  = \int_{\delta\mathcal{S}(\mathbf{x})} \nabla_{\mathbf{x}} \mathbf{z}^{T} \mathbf{n} \mathbf{F}ds
\end{equation}
where $\mathbf{n}$ is the unit vector outwards normal to the boundary and $ds$ is a infinitesimal contour of the boundary if $n_{\mathbf{z}} = 2$, while it is the infinitesimal area of the boundary surface if $n_{\mathbf{z}} = 3$.

\section{Some useful integrals}\label{sec:useful}
In this section, we will provide some useful identities of integrals which will be used to get \eqref{eq:cov_TG_nD} from \eqref{eq:cov_TG_nD_1}. The first one is
\begin{align}
&\int_{0}^{\pi} \sin^{n} \theta d \theta = \int_{0}^{\pi} \sin^{n-2} \left(1 - \cos^{2}\theta \right) d \theta\nonumber \\
& = \int_{0}^{\pi} \sin^{n-2} \theta d \theta - \int_{0}^{\pi} \frac{1}{n-1} \cos\theta d \sin^{n-1} \theta \nonumber \\
& = \int_{0}^{\pi}\! \!\sin^{n-2} \theta d \theta -\nonumber\\
&\qquad \frac{1}{n-1}\left(\sin^{n-1} \theta \cos\theta\Bigr\rvert_{\theta = 0}^{\pi}\! \!\!-\!\!  \int_{0}^{\pi}\!\! \!\sin^{n-1} \theta d \cos\theta \right)\nonumber\\
& = \int_{0}^{\pi} \sin^{n-2} \theta d \theta - \frac{1}{n-1} \sin^{n} \theta d\theta
\end{align}
Therefore,
\begin{align}\label{eq:integral_equality_1}
\int_{0}^{\pi} \sin^{n} \theta d \theta = \frac{n-1}{n} \int_{0}^{\pi} \sin^{n-2} \theta d\theta
\end{align}
The second is 
\begin{align}\label{eq:integral_equality_2}
&\int_{r = 0}^{a} \exp(-\frac{r^2}{2})r^{n+1} dr = - \int_{r = 0}^{a} r^{n} d\exp(-\frac{r^2}{2})\nonumber\\
= &\  r^n \exp(-\frac{r^2}{2})\Bigr\rvert_{r = 0}^{a} + n \int_{r = 0}^{a} \exp(-\frac{r^2}{2})r^{n-1} dr\nonumber\\
 = &\  -a^n\exp(-\frac{a^2}{2}) +  n \int_{r = 0}^{a} \exp(-\frac{r^2}{2})r^{n-1} dr
\end{align}
The third is
\begin{align}\label{eq:integral_equality_3}
\int_{\theta= 0}^{2\pi} \sin^2\theta d\theta = \pi
\end{align}

\bibliographystyle{IEEEtranS}
\bibliography{CRLLB}

\end{document}